\documentclass[a4paper,11pt,twoside]{amsart}
\usepackage[english]{babel}
\usepackage[utf8]{inputenc}

\usepackage[a4paper,inner=3cm,outer=3cm,top=4cm,bottom=4cm,pdftex]{geometry}
\usepackage{fancyhdr}
\usepackage{titlesec}
\usepackage{color}
\usepackage{bold-extra}
\usepackage{ mathrsfs }
\titleformat{\section}{\normalfont\scshape\centering}{\thesection}{1em}{}
  \titleformat{\subsection}{\bfseries}{\thesubsection}{1em}{}
  
\usepackage{comment}
\usepackage{graphics}
\usepackage{aliascnt}
\usepackage[pdftex,citecolor=green,linkcolor=red]{hyperref}

\usepackage{amsmath}
\usepackage{amsfonts}
\usepackage{amssymb}
\usepackage{amsthm}
\usepackage{comment}
\usepackage{mathtools}
\newcommand{\qedd}{\hfill \ensuremath{\Box}}
\newtheorem{theorem}{Theorem}[section]

\newtheorem{lemma}[theorem]{Lemma}

\theoremstyle{definition}
\newtheorem{definition}[theorem]{Definition}
\newtheorem{remark}[theorem]{Remark}
\newtheorem{conjecture}[theorem]{Conjecture}
\numberwithin{equation}{section}

\title{ON BINARY CORRELATIONS OF MULTIPLICATIVE FUNCTIONS}
\author{Joni Ter\"av\"ainen}
\date{}
\begin{document}

\begin{abstract} We study logarithmically averaged binary correlations of bounded multiplicative functions $g_1$ and $g_2$. A breakthrough on these correlations was made by Tao, who showed that the correlation average is negligibly small whenever $g_1$ or $g_2$ does not pretend to be any twisted Dirichlet character, in the sense of the pretentious distance for multiplicative functions. We consider a wider class of real-valued multiplicative functions $g_j$, namely those that are uniformly distributed in arithmetic progressions to fixed moduli. Under this assumption, we obtain a discorrelation estimate, showing that the correlation of $g_1$ and $g_2$ is asymptotic to the product of their mean values. We derive several applications, first showing that the numbers of large prime factors of $n$ and $n+1$ are independent of each other with respect to logarithmic density. Secondly, we prove a logarithmic version of the conjecture of Erd\H{o}s and Pomerance on two consecutive smooth numbers. Thirdly, we show that if $Q$ is cube-free and belongs to the Burgess regime $Q\leq x^{4-\varepsilon}$, the logarithmic average around $x$ of the real character $\chi \pmod{Q}$ over the values of a reducible quadratic polynomial is small.
\end{abstract}
\maketitle

\section{Introduction}\label{sec:intro}

Let $\mathbb{D}=\{z\in \mathbb{C}:\,\, |z|\leq 1\}$ be the unit disc of the complex plane, and let $g_1,g_2:\mathbb{N}\to \mathbb{D}$ be multiplicative functions. We consider the logarithmically averaged binary correlations
\begin{align}\label{eqtn1}
\frac{1}{\log x}\sum_{n\leq x}\frac{g_1(n)g_2(n+h)}{n},    
\end{align}
with $h\neq 0$ a fixed integer and $x$ tending to infinity. If $h<0$ in \eqref{eqtn1} , we can extend $g_1$ and $g_2$ arbitrarily to the negative integers, since this affects \eqref{eqtn1} only by $o(1)$.\\

In a recent breakthrough work,  Tao \cite{tao} showed that the correlation \eqref{eqtn1} is $o(1)$ as $x\to \infty$, provided that at least one of the two functions $g_j$ does not \emph{pretend} to be a twisted Dirichlet character, in the sense that
\begin{align}\label{eqtn2}
\liminf_{X\to \infty}\inf_{|t|\leq X}\mathbb{D}(g_j,\chi(n)n^{it};X)=\infty,
\end{align}
for all fixed Dirichlet characters $\chi$, with the pretentious distance $\mathbb{D}(\cdot)$ measured by
\begin{align}\label{eqq69}
\mathbb{D}(f,g;X):=\left(\sum_{p\leq X} \frac{1-\textnormal{Re}(f(p)\overline{g(p)})}{p}\right)^{\frac{1}{2}}.   
\end{align}
The main theorem in \cite{tao} that \eqref{eqtn1} is $o(1)$ under the non-pretentiousness assumption \eqref{eqtn2} is a logarithmically averaged version of the binary case of a conjecture of Elliott. Elliott's original conjecture \cite{elliott}, \cite{elliott2} (in the slightly corrected form presented in \cite{mrt}) states that for any integer $k\geq 1$, any multiplicative functions $g_1,\ldots, g_k:\mathbb{N}\to \mathbb{D}$ and any distinct integer shifts $h_1,\ldots, h_k$ we have the discorrelation estimate
\begin{align}\label{eqtn2a}
\frac{1}{x}\sum_{n\leq x}g_1(n+h_1)\cdots g_k(n+h_k)=o(1)    
\end{align}
as $x\to \infty$, provided that at least one of the $g_j$ satisfies the non-pretentiousness assumption
\eqref{eqtn2}.\footnote{In the case where the functions $g_j$ are allowed to depend on $x$, one needs a slightly stronger pretentiousness hypothesis; see \cite{tao_blog}.} The case $k=1$ is known as Halász's theorem \cite{halasz}. Already for $k=2$, there is not much progress towards the non-logarithmic version of Elliott's conjecture (see though \cite{elliott}). However, if one averages \eqref{eqtn2a} over the shifts $h_1,\ldots, h_k\in [1,H]$, with $H=H(x)$ tending to infinity with any speed, then Elliott's conjecture holds on average by the work of Matom\"aki, Radziwi\l{}\l{} and Tao \cite{mrt}. This was generalized by Frantzikinakis \cite{frantz-average} to averages along independent polynomials. In the case of logarithmically averaged correlations, there has been a lot of recent progress, initiated by \cite{tao}; see \cite{tao_sarnak}, \cite{frantz-chowla}, \cite{tao-teravainen}.\\

We study in this paper the same logarithmically averaged correlation \eqref{eqtn1} as Tao studied in \cite{tao}, but for a wider class of real-valued multiplicative functions (in \cite{tao} one works also with complex-valued functions). The multiplicative functions $g_j:\mathbb{N}\to [-1,1]$  that we consider are uniformly distributed in residue classes to fixed moduli. Many of the most interesting bounded multiplicative functions have such a uniform distribution property; in particular, the Liouville function $\lambda$ and the indicator function of $x^{a}$-smooth numbers up to $x$ have that property. Also the real primitive Dirichlet character $\chi_Q \pmod{Q}$ will be seen to be uniformly distributed in arithmetic progressions on $[x,2x]$, provided that the modulus $Q$ grows neither too slowly nor too rapidly in terms of $x$. Indeed, many of the applications of our main theorem concern consecutive smooth (or friable) numbers or quadratic residues.\\

The uniformity assumption we require of multiplicative functions is as follows. 

\begin{definition}[Uniformity assumption] \label{def1} Let $x\geq 1$, $1\leq Q\leq x$ and $\eta>0$. For a function $g:\mathbb{N}\to \mathbb{D}$, we write $g\in \mathcal{U}(x,Q,\eta)$ if we have the estimate
\begin{align*}
\bigg|\frac{1}{x}\sum_{\substack{x\leq n\leq 2x\\n\equiv a \pmod q}}g(n)-\frac{1}{qx}\sum_{x\leq n\leq 2x}g(n)\bigg|\leq\frac{\eta}{q}\quad \text{for all}\quad 1\leq a\leq q\leq Q.  
\end{align*}
\end{definition}

\begin{remark} Note that in this definition we do not send $x$ to infinity (but naturally we want $x$ to be large). The fact that Definition \ref{def1} is not an asymptotic relation is important, since later we shall to apply it to $g(n)=1_{n\leq x,\, n\,\textnormal{is}\,x^{a}\,\textnormal{-smooth}}$, which is a function dependent on $x$. 
\end{remark}

\begin{remark} \label{rmk2} Let $g:\mathbb{N}\to [-1,1]$ be a non-pretentious multiplicative function, in the sense that for some small $\varepsilon>0$ and some large $x$ we have
\begin{align*}
\inf_{|t|\leq x}\mathbb{D}(g,\chi(n)n^{it};x)\geq \varepsilon^{-10}    
\end{align*}
for all Dirichlet characters $\chi$ of modulus $\leq \varepsilon^{-10}$. By expressing the condition $n\equiv a \pmod{q}$ in Definition \ref{def1} in terms of Dirichlet characters (after reducing to a coprime residue class), and applying Halász's theorem, one sees that $g\in \mathcal{U}(x,\varepsilon^{-1},\varepsilon)$. Therefore, the collection of uniformly distributed real-valued multiplicative functions $g:\mathbb{N}\to [-1,1]$ contains all non-pretentious real functions.
\end{remark}

We will use the notation $o_{\varepsilon\to 0}(1)$ to denote a quantity depending on $\varepsilon$ and tending to $0$ as $\varepsilon\to 0$, uniformly with respect to all other parameters. With this notation, our main theorem asserts the following.

\begin{theorem}\label{theo_bincorr} Let a small real number $\varepsilon>0$, a fixed integer $h\neq 0$, and a function $\omega:\mathbb{R}_{\geq 1}\to \mathbb{R}$ with $1\leq \omega(X)\leq \log(3X)$ and $\omega(X) \xrightarrow{X\to \infty} \infty$ be given. Let $x\geq x_0(\varepsilon, h, \omega)$. Then, for any multiplicative functions $g_1, g_2:\mathbb{N}\to [-1,1]$  such that $g_1\in \mathcal{U}(x,\varepsilon^{-1},\varepsilon)$, we have
\begin{align*}
\frac{1}{\log \omega(x)}\sum_{\frac{x}{\omega(x)}\leq n\leq x}\frac{g_1(n)g_2(n+h)}{n}=\left(\frac{1}{x}\sum_{x\leq n\leq 2x}g_1(n)\right)\left(\frac{1}{x}\sum_{x\leq n\leq 2x}g_2(n)\right) +o_{\varepsilon\to 0}(1).
\end{align*}
\end{theorem}

\begin{remark}\label{rmk1}
Theorem \ref{theo_bincorr} can be viewed as stating that the functions $g_1$ and $g_2$ do not correlate with the shifts of each other. Note that in the case where the mean values of $g_1$ and $g_2$ on $[x,2x]$ are not $o(1)$, Theorem \ref{theo_bincorr} is not covered by the logarithmically averaged Elliott conjecture from \cite{tao}.   
\end{remark}

\begin{remark}
In the case where $g_1$ and $g_2$ are complex-valued, one does not always have the conclusion of Theorem \ref{theo_bincorr}. Namely, take $g_1(n)=n^{it}$ and $g_2(n)=n^{iu}$ for some $t,u\neq 0$ with $t+u\neq 0$. One easily sees that $g_1$ and $g_2$ are uniformly distributed in arithmetic progressions, and by partial summation the shifted product $g_1(n)g_2(n+1)=n^{i(t+u)}+o(1)$ has logarithmic mean value $o(1)$ on $[\frac{x}{\omega(x)},x]$. However, by the simple estimate $\frac{1}{x}\sum_{n\leq x}n^{it}=\frac{x^{it}}{1+it}+o(1)$, the product of the mean values of $g_1$ and $g_2$ on $[x,2x]$ is an oscillating function.
\end{remark}

\begin{remark}\label{rmk3} Although the statement of Theorem \ref{theo_bincorr} does not hold for all complex-valued multiplicative functions, one could show that it continues to hold if $g_1,g_2:\mathbb{N}\to \mathbb{D}$ take values in the roots of unity of fixed order. Indeed, the only places in the proof of the main theorem where real-valuedness plays a role are Lemmas \ref{le_stability}, \ref{le_C} and \ref{le_mr}. The first two lemmas could be proved also for functions $g_j$ taking values in the roots of unity of bounded order by applying a standard generalization of \cite[Lemma C.1]{mrt} to such functions. For Lemma \ref{le_mr}, one would also apply this generalization of \cite[Lemma C.1]{mrt} together with an extension of \cite[Theorem 3]{matomaki-radziwill} to multiplicative functions taking a bounded number of complex values. For this last extension, one notes that the only place in the proof of \cite[Proposition 1]{matomaki-radziwill} where real-valuedness is used is \cite[Lemma 3]{matomaki-radziwill}, and this lemma can also be made to work for functions taking values in the roots of unity of fixed order. We leave the details to the interested reader. 
\end{remark}

\begin{remark}
The bound $\omega(X)\leq \log (3X)$ in Theorem \ref{theo_bincorr} is not restrictive in reality, since if one wants an asymptotic formula for the logarithmic correlation over the interval $[1,x]$, say, one can sum together the asymptotics for the correlations over $[\frac{y}{\log(3y)},y]$ for various $y\leq x$. It is nevertheless necessary for technical reasons to have an upper bound on $\omega(X)$ in the main theorem, since otherwise the asymptotic would not be valid for example for the correlations of the indicator function of $x^{a}$-smooth numbers.
\end{remark}

\begin{remark} One could prove the same correlation bound for the more general logarithmic averages  of $g_1(a_1n+h_1)g_2(a_2n+h_2)$ with $(a_1,h_1)=(a_2,h_2)=1$ and $a_1,a_2\geq 1$ and $h_1,h_2$ fixed integers. This is due to the fact that the main theorem in \cite{tao} deals with such correlations. To avoid complicating the notations, however, we deal with the case $a_1=a_2=1$ here. 
\end{remark}

One might wonder at first why in the asymptotic formula in Theorem \ref{theo_bincorr} one side of the formula involves the values of the functions $g_j$ on $[\frac{x}{\omega(x)},x]$, whereas the other side only involves the values on $[x,2x]$. However, by a result we present in Appendix \ref{a:a}, essentially due to Granville and Soundararjan \cite{gs-spectrum}, the mean value of a real-valued multiplicative function is almost the same over the intervals $[\frac{x}{\omega(x)},x]$ and $[x,2x]$, explaining the phenomenon.\\ 

Owing to Remark \ref{rmk2}, the main theorem contains as a special case the logarithmically averaged binary Elliott conjecture from \cite{tao}. This is not surprising, since we use the same proof method. Of course, our interest lies in those cases where the functions $g_1$ and $g_2$ are pretentious (in the sense that \eqref{eqtn2} fails) but still satisfy our uniformity assumption.\\

It was recently shown by Klurman \cite{klurman} that one can obtain an asymptotic formula for the $k$-point correlations 
\begin{align*}
\frac{1}{x}\sum_{n\leq x}f_1(n+h_1)\cdots f_k(n+h_k)
\end{align*}
for any integers $h_1,\ldots h_k$, when $f_1,\ldots, f_k:\mathbb{N}\to \mathbb{D}$ are pretentious multiplicative functions, in the sense that $\mathbb{D}(f_j,\chi_j(n)n^{it_j};x)\ll 1$ for some characters $\chi_j$. This result does not imply Theorem \ref{theo_bincorr}, however, since our theorem is in a non-asymptotic form, allowing the multiplicative functions $g_1$ and $g_2$ to strongly depend on $x$. Indeed, allowing the multiplicative functions $g_j$ to depend on $x$ is crucial for applications to smooth numbers and to Burgess-type bounds. The asymptotic formula in \cite{klurman} is a sieve-theoretic product of local mean values, but one cannot express the density of smooth numbers as such a product.

\subsection{Applications of the main theorem}

We have a number of corollaries to Theorem \ref{theo_bincorr}. To state them, we recall the notion of logarithmic density of a set of integers.

\begin{definition} The \textit{logarithmic density} of a set $A\subset \mathbb{N}$ is
\begin{align*}
\delta(A)=\lim_{x\to \infty}\frac{1}{\log x}\sum_{\substack{n\leq x\\n\in A}}\frac{1}{n},    
\end{align*}
whenever it exists. 
\end{definition}

We will prove using Theorem \ref{theo_bincorr} the following theorem about the largest prime factors of consecutive integers.

\begin{theorem}[Independence of the number of large prime factors of $n$ and $n+1$] \label{theo_omega} Let $\omega_{> y}(n):=|\{p>y:\,\, p\mid n\}|$ be the number of prime factors of $n$ that are larger than $y$. Then, for any real numbers $a,b\in (0,1)$ and any integers $0\leq k< \frac{1}{a}$, $0\leq \ell< \frac{1}{b}$, we have
\begin{align*}
&\delta(\{n\in \mathbb{N}:\,\, \omega_{>n^{a}}(n)=k,\,\, \omega_{> n^{b}}(n+1)=\ell\})\\
&=\delta(\{n\in \mathbb{N}:\,\, \omega_{>n^{a}}(n)=k\})\cdot \delta(\{n\in \mathbb{N}:\,\, \omega_{> n^{b}}(n)=\ell\}).
\end{align*}
Moreover, under the same assumptions, the set $\{n\in \mathbb{N}:\,\, \omega_{>n^{a}}(n)=k,\,\, \omega_{> n^{b}}(n+1)=\ell\}$ has positive asymptotic lower density.
\end{theorem}

\begin{remark} From the proof of Theorem \ref{theo_omega} in Section \ref{sec:apps}, we can easily deduce a discorrelation estimate for the ''truncated Liouville function'' $\lambda_{>y}(n)$, which is a multiplicative function taking the value $+1$ at the primes $p\leq y$ and $-1$ at the primes $p>y$. This estimate takes the form 
\begin{align}\label{eqq110}
\frac{1}{\log x}\sum_{n\leq x}\frac{\lambda_{>x^{\varepsilon}}(n)\lambda_{>x^{\varepsilon}}(n+1)}{n}=o_{\varepsilon\to 0}(1).    
\end{align}
for  $\varepsilon\in (0,1)$ and $x\geq x_0(\varepsilon)$. This result may be compared with that of Daboussi and Sárk\H{o}zy \cite{daboussi-sarkozy} and Mangerel \cite{mangerel-chowla}, which states that if we define $\lambda_{<y}(n)$ as the completely multiplicative function taking the value $-1$ at the primes $p<y$ and $+1$ at the primes $p\geq y$ (so that $\lambda_{<y}(p)$ has the opposite sign as $\lambda_{>y}(p)$), then 
\begin{align}\label{eqq93}
\frac{1}{x}\sum_{n\leq x}\lambda_{<x^{\varepsilon}}(n)\lambda_{<x^{\varepsilon}}(n+1)=o_{\varepsilon\to 0}(1); 
\end{align}
moreover, they proved this in a quantitative form. The proof of \eqref{eqq93} is based on sieve theory and is very different from the proof of \eqref{eqq110}.
\end{remark}

Our next applications concern smooth numbers, so we introduce the function $P^{+}(n)$, whose value is the largest prime factor of the positive integer $n\geq 2$ (and  $P^{+}(1)=1$). We say that a number $n$ is $y$\emph{-smooth} if $P^{+}(n)\leq y$. The simultaneous distribution of the function $P^{+}(\cdot)$ at consecutive integers is the subject of several conjectures. There is for instance a conjecture of Erd\H{o}s and Pomerance \cite{erdos-pomerance}, asserting that the largest prime factors of $n$ and $n+1$ are independent events.

\begin{conjecture}[Erd\H{o}s--Pomerance] \label{conj_pomerance}
For any $a,b\in (0,1)$, the asymptotic density of the set
\begin{align}\label{eqtn5}
\{n\in \mathbb{N}:\,\, P^{+}(n)\leq n^{a},\,\, P^{+}(n+1)\leq n^{b}\}    \end{align}
exists and equals $\rho(\frac{1}{a})\rho(\frac{1}{b})$, where $\rho(\cdot)$ is the \emph{Dickmann function} (see \cite[Section 1]{ht}).
\end{conjecture}

What we are able to prove, taking $k=\ell=0$ in Theorem \ref{theo_omega}, is a logarithmic version of the conjecture.

\begin{theorem}\label{theo_density} Conjecture \ref{conj_pomerance} holds when asymptotic density is replaced with logarithmic density; that is, for any $a,b\in (0,1)$ we have
\begin{align*}
\delta(\{n\in \mathbb{N}:\,\, P^{+}(n)\leq n^{a},\,\, P^{+}(n+1)\leq n^{b}\})=\rho\left(\frac{1}{a}\right)\rho\left(\frac{1}{b}\right).    
\end{align*}
\end{theorem}

A closely related conjecture,  formulated in the correspondence of Erd\H{o}s and Turán in the 1930s (see \cite[pp. 100-101]{erdos-turan}, \cite{erdos-conj}, \cite[Section 1]{pollack}) is that the distribution of $(P^{+}(n),P^{+}(n+1))$ is symmetric.

\begin{conjecture}[Erd\H{o}s-Turán] \label{conj_turan} The asymptotic density of the set
\begin{align}\label{eqtn4}
\{n\in \mathbb{N}:\,\, P^{+}(n)<P^{+}(n+1)\}    
\end{align}
exists and equals $\frac{1}{2}$.
\end{conjecture}

There has been some progress towards this conjecture.  Erd\H{o}s and Pomerance \cite{erdos-pomerance} showed that the lower asymptotic density of the set in \eqref{eqtn4} is positive (in fact, at least $0.0099$). The lower bound for the density was improved to $0.05544$ by de la Bretèche, Pomerance and Tenenbaum \cite{bpt}, to $0.1063$ by Wang \cite{wang}, and a further improvement to $0.1356$ was given by Wang in \cite{wang2}.\\

We can prove Conjecture \ref{conj_turan} if asymptotic density is again replaced with logarithmic density.

\begin{theorem}\label{theo_erdos} Conjecture \ref{conj_turan} holds when asymptotic density is replaced with logarithmic density; that is,
\begin{align*}
\delta(\{n\in \mathbb{N}:\,\, P^{+}(n)<P^{+}(n+1)\})=\frac{1}{2}. \end{align*}
\end{theorem}

In fact, Theorem \ref{theo_density} implies Theorem \ref{theo_erdos}, via the following theorem, which was also conjectured by Erd\H{o}s \cite{erdos-conj} in the case of asymptotic density.\footnote{Erd\H{o}s conjectured the existence of the density of integers $n$ for which $P^{+}(n+1)>P^{+}(n)\cdot n^{\alpha}$.}

\begin{theorem}\label{theo_integral} Let $\alpha \in [0,1]$ be a real number. Let $u(x):=\frac{\rho(\frac{1}{x}-1)}{x}$ for $x\in (0,1)$, where $\rho$ is the Dickmann function. Then we have
\begin{align}\label{eqq53}
\delta(\{n\in \mathbb{N}:\,\, P^{+}(n+1)>P^{+}(n)\cdot n^{\alpha}\})=\int_{T_{\alpha}}u(x)u(y)\, dx\, dy,
\end{align}
where $T_{\alpha}$ is the triangular domain $\{(x,y)\in [0,1]^2:\,\, y\geq x+\alpha\}$. In particular, the logarithmic density above exists.
\end{theorem}

\begin{remark} The appearance of the function $u(\cdot)$ is to be expected in Theorem \ref{theo_integral}, since $u$ is the derivative of $x\mapsto \rho(\frac{1}{x})$, with the latter function expressing the probability that $P^{+}(n)\leq n^{x}$. 
\end{remark}

We will prove Theorem \ref{theo_omega}, and consequently Theorem \ref{theo_density}, in Section \ref{sec:apps}, where we will also see that Theorem \ref{theo_integral} quickly follows from the latter theorem. With Theorem \ref{theo_integral} available, Theorem \ref{theo_erdos} follows by taking $\alpha=0$ and noting that then the integral in \eqref{eqq53} is symmetric in $x$ and $y$, implying that its value is $\frac{1}{2}$. For the details, see Section \ref{sec:apps}.\\

We can also prove another approximation to Conjecture \ref{conj_pomerance}. This was obtained earlier by Hildebrand \cite{hildebrand}, using a combinatorial method, in the special case $(a,b)=(c,d)$ (Hildebrand's proof also applies to so-called stable sets, with power-smooth numbers being an example of such a set).  The following theorem also implies a result of Wang \cite[Théorème 2]{wang2} on the integers $n\leq x$ with $P_{y}^{+}(n)<P_{y}^{+}(n+1)$ having a positive density, where $P_y^{+}(n)=\max\{p\leq y:\,p\mid n\}$ and $y\geq x^{\varepsilon}$.

\begin{theorem}\label{theo_hildebrand} Let $a,b,c,d\in (0,1)$ be real numbers with $a<b$ and $c<d$. Then the set
\begin{align*}
\{n\in \mathbb{N}:\,\, n^{a}\leq P^{+}(n)\leq n^{b},\,\, n^{c}\leq P^{+}(n+1)\leq n^{d}\}    
\end{align*}
has positive asymptotic lower density.
\end{theorem}

 Note that Theorem \ref{theo_hildebrand} is not implied by Theorem \ref{theo_density}, as there are sets of positive logarithmic density having zero asymptotic lower density. Nevertheless, the proof we use for the latter theorem also works for the former, owing to the presence of an arbitrarily slowly growing function $\omega(X)$ in Theorem \ref{theo_bincorr}.\\

Since we can prove satisfactory results for the distribution of the largest prime factor function $P^{+}(\cdot)$ at two consecutive integers, it is natural to ask about the distribution of $P^{+}(\cdot)$ also at longer strings of consecutive integers. A conjecture of De Koninck and Doyon \cite{kd} states the following.

\begin{conjecture}[De Koninck and Doyon] \label{conj_konick}Let $k\geq 2$ be an integer and $(a_1,\ldots a_k)$ any permutation of the set $\{1,2,\ldots, k\}$. Then the set 
\begin{align}\label{eqtn30}
\{n\in \mathbb{N}:\,\, P^{+}(n+a_1)<\ldots<P^{+}(n+a_k)\}    
\end{align}
has an asymptotic density, and it equals $\frac{1}{k!}$. 
\end{conjecture}

The case $k=2$ of this is the earlier mentioned Conjecture \ref{conj_turan} of Erd\H{o}s and Turán. Little is known about this conjecture for $k\geq 3$; it is not even known that the sets in \eqref{eqtn30} have positive asymptotic lower density. Recently, Wang \cite{wang2} proved a result about orderings of $P^{+}(\cdot)$ at consecutive integers, showing that
\begin{align}\label{eqq42}
 P^{+}(n+i)<\min_{\substack{j\leq J\\j\neq i}} P^{+}(n+j)\quad \textnormal{and}\quad  P^{+}(n+i)>\max_{\substack{j\leq J\\j\neq i}} P^{+}(n+j) 
\end{align}
hold with positive asymptotic lower density for any $J\geq 3$ and $1\leq i\leq J$. The method of \cite{wang2} is based on the linear sieve and Bombieri-Vinogradov type estimates for smooth numbers. Applying Theorem \ref{theo_bincorr} together with the Matom\"aki-Radziwi\l{}\l{} theorem \cite{matomaki-radziwill} on multiplicative functions in short intervals (and using the method of \cite{matomaki-sign}), we can give a different proof of the $J=3$ case of Wang's result. We leave the details of this special case of \eqref{eqq42} to the interested reader.\\

As our last application, we study character sums along the values of a reducible quadratic polynomial $n(n+h)$. A famous result of Burgess \cite{burgess} states that for any non-principal Dirichlet character $\chi$ modulo $Q$ we have
\begin{align*}
\sum_{y\leq n\leq y+x}\chi(n)\ll_{r,\varepsilon} x^{1-\frac{1}{r}}Q^{\frac{r+1}{4r^2}+\varepsilon},  
\end{align*}
whenever $r\in \mathbb{N}$ and $Q$ is cube-free (that is, $p^3\nmid Q$ for all primes $p$). In particular, we have the important special case
\begin{align}\label{eqq70}
\sum_{n\leq x}\chi_{Q}(n)=o(x), \quad 3\leq Q\leq x^{4-\varepsilon}  
\end{align}
for cube-free values of $Q$, where $\chi_{Q}$ is a real primitive Dirichlet character modulo $Q$. Using Theorem \ref{theo_bincorr}, we can prove that a variant of the estimate \eqref{eqq70} continues to hold for character sums over the values of a reducible quadratic polynomial.

\begin{theorem}[Character sums over $n(n+h)$ in the Burgess regime]\label{theo_burgess} Let $\varepsilon>0$ be small, $h\neq 0$ a fixed integer, and $1\leq \omega(X)\leq \log(3X)$ any function tending to infinity. For $x\geq x_0(\varepsilon,h,\omega)$, let $Q=Q(x)\leq x^{4-\varepsilon}$ be a cube-free natural number  with $Q(x)\xrightarrow{x\to \infty}\infty$. Then, the real primitive Dirichlet character $\chi_Q$ modulo $Q$ satisfies the estimate
\begin{align*}
\frac{1}{\log \omega(x)}\sum_{\frac{x}{\omega(x)}\leq n\leq x} \frac{\chi_Q(n(n+h))}{n}=o(1). \end{align*}
Moreover, if $Q$ is as before and QNR stands for quadratic nonresidue\footnote{We say that $n$ is a quadratic nonresidue $\pmod{Q}$ if $\chi_{Q}(n)=-1$.}, we have
\begin{align}\label{eqq73}
\frac{1}{\log x}\sum_{\substack{n\leq x\\n,\,n+1\, \textnormal{QNR}\,\, \hspace{-0.2cm}\pmod{Q}}}\frac{1}{n}=\frac{1}{4}\prod_{p\mid Q}\left(1-\frac{2}{p}\right)+o(1) \end{align}
and
\begin{align}\label{eqq74}
\frac{1}{x}|\{n\leq x:\,\, n\,\, \textnormal{and}\,\, n+1\,\, \textnormal{QNR}\,\, \hspace{-0.2cm}\pmod{Q}\}|\gg \prod_{p\mid Q}\left(1-\frac{2}{p}\right) .   
\end{align}
\end{theorem}

\begin{remark} In light of Remark \ref{rmk3}, we could also prove Theorem \ref{theo_burgess} for primitive characters $\chi$ modulo $Q$ whose order is bounded (that is, characters $\chi$ such that $\chi^k$ is principal for some $k\ll 1$).
\end{remark}

This theorem is related to \cite[Problem 11]{shparlinski}, although there one asks for cancellation in the ordinary average instead of the logarithmic one, and one wants to take a maximum over $h\leq Q$ (but there $Q$ is restricted to primes and $Q\leq x^{2+\delta}$ for some small $\delta>0$). We also remark that in the much smaller range $Q=o( x^{2}/(\log x))$ and with $Q$ prime,  one can use the Weil bound \cite[Theorem 11.23]{iwaniec-kowalski} to prove the above estimate. In the same range $Q\leq x^{4-\varepsilon}$ as in Theorem \ref{theo_burgess}, it was shown by Burgess \cite{burgess2} that
\begin{align*}
x-\sum_{y\leq n\leq y+x}\chi_Q(n(n+h))\gg_{\varepsilon,h} x^{\frac{\varepsilon}{2}},    
\end{align*}
and the same estimate holds with $n(n+h)$ replaced by any polynomial that factorizes into linear factors and is not the square of another polynomial.\\

We note that Theorem \ref{theo_burgess} does not directly follow from the logarithmically averaged binary Elliott conjecture proved in \cite{tao}, since if the Vinogradov quadratic nonresidue conjecture\footnote{Namely, that for any $q\geq q(\varepsilon)$, there is a quadratic nonresidue $\pmod{q}$ on the interval $[1,q^{\varepsilon}]$.} failed, it would be the case that
\begin{align}\label{eqq71}
\mathbb{D}(\chi_Q,1;x)\ll 1.    
\end{align}
We of course do not expect \eqref{eqq71} to hold, but it cannot be ruled out with current knowledge. Furthermore, the correlation asymptotic in \cite{klurman} does not apply either to Theorem \ref{theo_burgess}, since the function $\chi_Q$ depends heavily on the length $x$ of the sum. Nevertheless, the function $\chi_Q$ has mean value $o(1)$ by the Burgess bound, and by a slight generalization of that, it also has mean $o(1)$ in fixed arithmetic progressions, which is what is required to apply Theorem \ref{theo_bincorr}. For the details of the proof of Theorem \ref{theo_burgess}, see Section \ref{sec:apps}.

\subsection{Structure of the paper}

The main theorem, Theorem \ref{theo_bincorr} will be proved in Sections \ref{sec:entropy} and \ref{sec:circle}. In the former of these sections, the entropy decrement argument from \cite{tao}, \cite{tao-teravainen} is deployed to replace the correlation average with a simpler, bilinear average. The proof of one lemma in Section \ref{sec:entropy}, concerning stability of mean values of multiplicative functions, is postponed to Appendix \ref{a:a}. In Section \ref{sec:circle}, we use circle method estimates and a short exponential sum estimate for multiplicative functions to show that the bilinear average we mentioned has the anticipated asymptotic formula, concluding the proof. The proof of this exponential sum estimate, which is a slight modification of the one by Matom\"aki, Radziwi\l{}\l{} and Tao \cite{mrt}, is left to Appendix \ref{a:b}. In Section \ref{sec:apps}, we apply Theorem \ref{theo_bincorr} to deduce the applications mentioned in the Introduction. Theorem \ref{theo_omega} will be proved first, and then Theorems \ref{theo_density} and \ref{theo_hildebrand} will be deduced from this. Theorems \ref{theo_integral} and \ref{theo_erdos} will in turn  follow from Theorem \ref{theo_density}. Theorem \ref{theo_burgess} will be deduced from the main theorem and the Burgess bound.

\subsection{Notation} The functions $g_1,g_2:\mathbb{N}\to [-1,1]$ are always multiplicative functions. The pretentious distance $\mathbb{D}(f,g;x)$ between two multiplicative functions is given by \eqref{eqq69}. We denote by $\mu(n)$ the M\"obius function, by $\varphi(n)$ the Euler totient function, and by $P^{+}(n)$ the largest prime factor of $n$, with the convention that $P^{+}(1)=1$. By $(a,b)$, we denote the greatest common divisor of $a$ and $b$. For a proposition $P(n)$, the indicator $1_{P(n)}$ is defined as $1$ if $P(n)$ is true and as $0$ if $P(n)$ is false. By $\delta(S)$ we denote the logarithmic density of $S\subset \mathbb{N}$, not to be confused with $\delta_1$, $\delta_2\in [-1,1]$, which are the mean values of $g_1$ and $g_2$, as defined in formula \eqref{eqq82}.\\

The variables $p,p_1,p_2,\ldots$ will always be primes. We reserve various letters, such as $d,k,\ell,m,n,q$ for positive integer quantities. The variables $x,y$ in turn will be understood to be large, whereas $\varepsilon>0$ will tend to zero. The integer $h\neq 0$ is always fixed, and the function $\omega:\mathbb{R}_{\geq 1}\to \mathbb{R}$ is a growth function satisfying $1\leq \omega(X)\leq \log(3X)$ and tending to infinity with $X$.\\

We use the standard Landau and Vinogradov asymptotic notations $O(\cdot), o(\cdot),\gg, \ll$, with the convention that the implied constants are absolute unless otherwise indicated. Thus for instance $o_{\varepsilon\to 0}(1)$ denotes a quantity depending on $\varepsilon$ and tending to $0$ as $\varepsilon\to 0$, uniformly with respect to all other involved parameters. All the logarithms in the paper will be to base $e$, and the function $\log_j x$ is the $j$th iterate of the logarithm function. The function $\exp_j x$ is analogously the $j$th iterate of $x\mapsto e^{x}$.

\subsection{Acknowledgments}

The author is grateful to Kaisa Matom\"aki and Terence Tao for various useful comments and discussions. He thanks the referee for careful reading of the paper and for useful comments. The author also thanks Zhiwei Wang for showing his preprint \cite{wang2} to the author. Thanks go also to Alexander Mangerel for a discussion on character sums. Part of this work was done while the author was visiting the Mathematical Sciences Research Institute (funded by NSF grant DMS-1440140) in spring 2017, and he thanks the institute for a stimulating atmosphere. The author was funded by UTUGS Graduate School and project number 293876 of the Academy of Finland.

\section{The entropy decrement argument and some reductions}\label{sec:entropy}

Given a function $\omega(X)$ having the same properties as in Theorem \ref{theo_bincorr}, we define for $a\in \mathbb{Z}$ the correlation sequence
\begin{align*}
f_{x,\omega}(a):=\frac{1}{\log \omega(x)}\sum_{\frac{x}{\omega(x)}\leq n\leq x}\frac{g_1(n)g_2(n+a)}{n}.    
\end{align*}
As was noted in the Introduction, one can define this equally well for $a<0$. Our task is then to show that if 
\begin{align}\label{eqq82}
\delta_1:=\frac{1}{x}\sum_{x\leq n\leq 2x}g_1(n),\quad \delta_2:=\frac{1}{x}\sum_{x\leq n\leq 2x}g_2(n)    
\end{align}
are the mean values of $g_1$ and $g_2$ (which depend on $x$), then $|f_{x,\omega}(h)-\delta_1\delta_2|=o_{\varepsilon\to 0}(1)$ under the assumptions of Theorem \ref{theo_bincorr}. By replacing $\varepsilon$ with $1/\exp_2(\varepsilon^{-2})$ in Theorem \ref{theo_bincorr}, with $\exp_2$ the second iterated exponential, we may in fact assume that
\begin{align}
g_1\in \mathcal{U}(x,\exp_2(\varepsilon^{-2}),1/\exp_2(\varepsilon^{-2}));    
\end{align}
we do this for notational convenience. We may also assume that  $|h|\leq \varepsilon^{-1}$, since $h$ is fixed in Theorem \ref{theo_bincorr} and  $\varepsilon$ is small.\\

We will average $f_{x,\omega}(h)$ over the primes belonging to a small scale using multiplicativity, and then apply the entropy decrement argument to relate $f_{x,\omega}(h)$ to a bilinear analogue $\frac{\log P}{P}\sum_{p\sim P} \frac{g_1(p)^{-1}g_2(p)^{-1}}{p} f_{x,\omega}(ph)$ of the same sum (this is the same approach as in Tao's paper \cite{tao}, and in the later works \cite{tao_sarnak}, \cite{tao-teravainen}). As in \cite{tao}, we will then apply the circle method and establish a slight variant of the short exponential sum estimate for multiplicative functions, due to Matom\"aki, Radziwi\l{}\l{} and Tao \cite{mrt}, to finish the proof. Since Theorem \ref{theo_bincorr} involves both pretentious and non-pretentious functions $g_j$, we need to make a distinction between them in certain parts of the argument. We will also separate the case where $|g_1(p)g_2(p)|$ is small for many primes $p$ from the opposite case, since expressions such as $g_1(p)^{-1}g_2(p)^{-1}$ naturally appear in the proof. To deal with these distinctions for $g_j$, we will need the fact that the entropy argument works not only in infinitely many dyadic scales $[2^m,2^{m+1}]$, but in fact in almost all of them with respect to some measure. Such a strengthening was presented in \cite{tao-teravainen}. We begin with this entropy decrement argument.

\begin{lemma}[Entropy decrement argument]\label{le_entropy} Let $\varepsilon>0$ be small, $|h|\leq \varepsilon^{-1}$ an integer, $x\geq x_0(\varepsilon, h, \omega)$, and $\omega:\mathbb{R}_{\geq 1}\to \mathbb{R}$ a function with $1\leq \omega(X)\leq X$ and $\omega(X) \xrightarrow{X\to \infty} \infty$. Let $g_1,g_2:\mathbb{N}\to \mathbb{D}$ be $1$-bounded multiplicative functions and $c_p\in \mathbb{D}$ any complex numbers. Then for all $m\in \mathcal{M}\cap [1,\log_2 \omega(x)]$ we have
\begin{align*}
\frac{m\log 2}{2^m}\sum_{2^m\leq p<2^{m+1}}c_p g_1(p)g_2(p)\cdot f_{x,\omega}(h)=\frac{m\log 2}{2^m}\sum_{2^m\leq p<2^{m+1}}c_p f_{x,\omega}(ph)+o_{\varepsilon\to 0}(1),    
\end{align*}
with the set $\mathcal{M}\subset \mathbb{N}$ being independent of $c_p$ and being large in the sense that
\begin{align}\label{eqq10}
\sum_{\substack{m\geq 1\\m\not \in \mathcal{M}}}\frac{1}{m}\ll \varepsilon^{-10}.    
\end{align}
\end{lemma}

\textbf{Proof.} This follows from the proof of \cite[Theorem 3.6]{tao-teravainen}, but since that argument uses generalized limit functionals, we outline how it goes through without them. We also remark that without the density bound \eqref{eqq10} Lemma \ref{le_entropy} follows  from \cite[Section 3]{tao}, and that in \cite[Theorem 3.1]{tt-chowla} the lemma was proved in the special case of the Liouville function.\\

We may assume that $m\geq \varepsilon^{-1}$ for all $m\in \mathcal{M}$, since removing the numbers $m<\varepsilon^{-1}$ from $\mathcal{M}$ alters the sum in \eqref{eqq10} by $\sum_{m<\varepsilon^{-1}}\frac{1}{m}\ll \varepsilon^{-1}$. We have the multiplicativity property $g_j(p)g_j(n)=g_j(pn)+O(1_{p\mid n})$ for any prime $p$, so for $2^m\leq p<2^{m+1}$ with $\varepsilon^{-1}\leq m\leq \log_2 \omega(x)$ we have
\begin{align*}
 g_1(p)g_2(p)\cdot f_{x,\omega}(h)&=\frac{1}{\log \omega(x)}\sum_{\frac{x}{\omega(x)}\leq n\leq x}\frac{g_1(pn)g_2(pn+ph)}{n}+O\bigg(\frac{1}{\log \omega(x)}\sum_{\substack{\frac{x}{\omega(x)}\leq n\leq x\\p\mid n(n+h)}}\frac{1}{n}\bigg)\\
&=\frac{1}{\log \omega(x)}\sum_{\frac{x}{\omega(x)}\leq n\leq x}\frac{g_1(pn)g_2(pn+ph)}{n}+O(\varepsilon)\\
&=\frac{1}{\log \omega(x)}\sum_{\frac{px}{\omega(x)}\leq n\leq px}\frac{g_1(n)g_2(n+ph)}{n}p1_{p\mid n}+O(\varepsilon)\\
&=\frac{1}{\log \omega(x)}\sum_{\frac{x}{\omega(x)}\leq n\leq x}\frac{g_1(n)g_2(n+ph)}{n}p1_{p\mid n}+O(\varepsilon),
\end{align*}
where the last step comes from estimating the terms $n\in [\frac{x}{\omega(x)},\frac{px}{\omega(x)}]$ and $n\in [x,px]$ trivially.\footnote{This is the part of the argument where it is crucial to work with logarithmic averaging.}\\

Define the modified functions $g_j^{(\varepsilon)}(n)$ by rounding $g_j(n)$ to the nearest element of the Gaussian lattice $\varepsilon \mathbb{Z}[i]$. Then, averaging over $p$ the above formula for $g_1(p)g_2(p)\cdot f_{x,\omega}(h)$, we get 
\begin{align}\label{eqq13}\begin{split}
&\frac{m\log 2}{2^m}\sum_{2^m\leq p<2^{m+1}}c_p g_1(p)g_2(p)\cdot f_{x,\omega}(h)\\
&=\frac{m\log 2}{2^m\log \omega(x)}\sum_{2^m\leq p<2^{m+1}}c_p\sum_{\frac{x}{\omega(x)}\leq n\leq x}\frac{g_1^{(\varepsilon)}(n)g_2^{(\varepsilon)}(n+ph)}{n}p1_{p\mid n}+O(\varepsilon).    
\end{split}
\end{align}
The concentration of measure argument in \cite{tao-teravainen} tells that we may replace $p1_{p\mid n}$ with $1+O(\varepsilon)$ in \eqref{eqq13}, provided that the random variables
\begin{align*}
\mathbf{X}_m:&=(g_r^{(\varepsilon)}(\mathbf{n}+j))_{1\leq r\leq 2,\,0\leq j\leq (1+|h|)2^{m+2}},\\  \mathbf{Y}_m:&=(\mathbf{n}\hspace{-0.2cm}\pmod p)_{2^m\leq p<2^{m+1}},\,\, \mathbf{Y}_{<m}:=(\mathbf{Y}_{m'})_{m'<m}  \end{align*}
enjoy the conditional mutual information\footnote{For the definition of conditional mutual information, see \cite[Section 2]{tao-teravainen}.} bound
\begin{align}\label{eqq14}
\mathbf{I}(\mathbf{X}_m:\mathbf{Y}_m|\mathbf{Y}_{<m})\leq \varepsilon^{4} \cdot \frac{2^m}{m}.
\end{align}
We thus need to show that the set $\mathcal{M}$ of $m$ for which \eqref{eqq14} holds satisfies \eqref{eqq10}. But this was shown in \cite[Proposition 3.5]{tao-teravainen} (see also Remark 3.7 there), so we obtain the claim.\qedd\\

Before utilizing Lemma \ref{le_entropy}, we will show that the quantities $\delta_1$ and $\delta_2$ in \eqref{eqq82} are the mean values of $g_1$ and $g_2$ also on many other intervals than $[x,2x]$. For this we use a slight generalization of a lemma due to Elliott \cite{elliott-lipschitz} and Granville and Soundarajan \cite[Proposition 4.1]{gs-spectrum}. Such results are also proved in Matthiesen's work \cite{matthiesen-fourier} in a more general setting.

\begin{lemma}[Stability of mean values of multiplicative functions]\label{le_stability} Let $g:\mathbb{N}\to [-1,1]$ be a real-valued multiplicative function, $x\geq 10$, and  $y\in [1,\log^{10} x]$ arbitrary. Then, for $a,q\in \mathbb{N}$ we have
\begin{align*}
\bigg|\frac{1}{x}\sum_{\substack{x\leq n\leq 2x\\n\equiv a\pmod{q}}}g(n)-\frac{1}{x/y}\sum_{\substack{\frac{x}{y}\leq n\leq \frac{2x}{y}\\n\equiv a \pmod{q}}}g(n)\bigg|\ll_q (\log x)^{-\frac{1}{400}}.    
\end{align*}
\end{lemma}

\textbf{Proof.} We prove this in Appendix \ref{a:a}.\qedd\\

Owing to the above lemma, we can show that the uniformity assumption on $g_1$ implies the seemingly stronger assumption that $g_1$ be uniformly distributed also on intervals $[\frac{x}{\omega(x)},x]$. For this purpose, we need the following definition.

\begin{definition}[Stronger uniformity assumption]\label{def2}  Let $1\leq Q\leq x$, $\eta>0$,  and $\delta \in \mathbb{C}$. Let  $\omega:\mathbb{R}_{\geq 1}\to \mathbb{R}$ be a function  with $1\leq \omega(X)\leq X$ for all $X\geq 1$. For a function $g:\mathbb{N}\to \mathbb{C}$, we write $g\in \mathcal{U}_{\omega}(x,Q,\eta,\delta)$ if we have the estimate
\begin{align*}
\bigg|\frac{1}{y}\sum_{\substack{y\leq n\leq 2y\\n\equiv a \pmod q}}g(n)-\frac{\delta}{q}\bigg|\leq\frac{\eta}{q}\quad \text{for all}\quad 1\leq a\leq q\leq Q\quad \textnormal{and}\quad \frac{x}{\omega(x)}\leq y\leq x.  \end{align*}
\end{definition}

With the above notation, if $\delta_1$ and $\delta_2$ are as in \eqref{eqq82} and $1\leq \omega(X)\leq \log(3X)$ is as in Theorem \ref{theo_bincorr}, by Lemma \ref{le_stability} we have
\begin{align}\label{eqq83}\begin{split}
g_1&\in \mathcal{U}_{\omega}(x,\exp_2(\varepsilon^{-2}),2/\exp_2(\varepsilon^{-2}),\delta_1),\\
g_2&\in \mathcal{U}_{\omega}(x,1,2/\exp_2(\varepsilon^{-2}),\delta_2)\,\, \textnormal{for}\,\, \omega(X)\leq \log^{10} X.
\end{split}
\end{align}
This property will be used several times in the rest of the proof of the main theorem. In particular, we have for all $y\in [x(\log x)^{-10}, x]$ the estimate
\begin{align*}
\sum_{y\leq n\leq 2y}g_j(n)=(\delta_j+O(\varepsilon))y,
\end{align*}
where, as always, the $O(\cdot)$ constant is absolute. Summing this over the dyadic intervals $[\frac{y}{2^{j+1}},\frac{y}{2^j}]$ for $j\geq 0$ and assuming that $y\geq x(\log(3x))^{-1}$, say, we get
\begin{align*}
\sum_{n\leq y}g_j(n)=(\delta_j+O(\varepsilon))y.    
\end{align*}
Subtracting this formula for two different lengths of summation, we see that
\begin{align*}
\sum_{y\leq n\leq z}g_j(n)=\delta_j(z-y)+O(\varepsilon z)    
\end{align*}
for all $x(\log(3x))^{-1}\leq y\leq z\leq 2x$. From this and partial summation, we obtain
\begin{align}\label{eqq84}
\frac{1}{\log \omega(x)}\sum_{\frac{x}{\omega(x)}\leq n\leq x}\frac{g_j(n)}{n}=\delta_j+O(\varepsilon)    
\end{align}
for $1\leq \omega(X)\leq \log(3X)$, which also will be utilized in what follows.\\

We return to applying the entropy argument. Defining the normalized correlation sequence
\begin{align*}
\tilde{f}_{x,\omega}(a):=\frac{1}{\log \omega(x)}\sum_{\frac{x}{\omega(x)}\leq n\leq x}\frac{(g_1(n)-\delta_1)(g_2(n+a)-\delta_2)}{n}
\end{align*}
and using the simple identity $XY=\delta_1\delta_2+\delta_1(Y-\delta_2)+\delta_2(X-\delta_1)+(X-\delta_1)(Y-\delta_2)$, we deduce from Lemma \ref{le_entropy} that
\begin{align}\label{eqq43}\begin{split}
&\frac{m\log 2}{2^m}\sum_{2^m\leq p<2^{m+1}}c_p g_1(p)g_2(p)\cdot f_{x,\omega}(h)\\
&=\delta_1\delta_2\frac{m\log 2}{2^m}\sum_{2^m\leq p<2^{m+1}}c_p +\frac{m\log 2}{2^m}\sum_{2^m\leq p<2^{m+1}}c_p \tilde{f}_{x,\omega}(ph)\\
&+O\bigg(\max_{r\in \{1,2\}}\frac{1}{\log \omega(x)}\bigg|\sum_{\frac{x}{\omega(x)}\leq n\leq x}\frac{g_r(n)-\delta_r}{n}\bigg|\bigg)+o_{\varepsilon \to 0}(1),\,\, m\in \mathcal{M}\cap [1,\log_2 \omega(x)].
\end{split}
\end{align}

Formula \eqref{eqq84} tells that the $O(\cdot)$ error term in \eqref{eqq43} is $o_{\varepsilon\to 0}(1)$. Then \eqref{eqq43} takes the form
\begin{align}\label{eqq5}\begin{split}
&\frac{m\log 2}{2^m}\sum_{2^m\leq p<2^{m+1}}c_p g_1(p)g_2(p)\cdot f_{x,\omega}(h)\\
&=\delta_1\delta_2\frac{m\log 2}{2^m}\sum_{2^m\leq p<2^{m+1}}c_p +\frac{m\log 2}{2^m}\sum_{2^m\leq p<2^{m+1}}c_p \tilde{f}_{x,\omega}(ph)+o_{\varepsilon \to 0}(1)\quad 
\end{split}
\end{align}
for $m\in \mathcal{M}\cap [1,\log_2 \omega(x)]$.\\

It is natural to predict that the average of the normalized correlation $\tilde{f}_{x,\omega}(h)$ in \eqref{eqq5} is small, and this is indeed what we will prove in Section \ref{sec:circle}. Before we deal with that term, we will consider the main term arising in \eqref{eqq5}. One would like to choose $c_p=g_1(p)^{-1}g_2(p)^{-1}$ there, since then the main term becomes just $\delta_1\delta_2+o_{\varepsilon\to 0}(1)$. However, it may be that $|g_j(p)|$ takes very small values (or even $0$), in which case $c_p$ would be unbounded. To avoid this, we prove two lemmas, the first of which tells that if the correlation average in Theorem \ref{theo_bincorr} is not negligibly small, then $|g_1(p)g_2(p)|\geq \frac{1}{2}$ for most primes $p$. The second lemma in turn tells that if $|\delta_1|$ and $|\delta_2|$ are not negligibly small, then $g_1(p)g_2(p)$ behaves like $1$ in most scales.

\begin{lemma}[Dealing with small values of $g_j(p)$]\label{le_A} Let the notations be as in Theorem \ref{theo_bincorr}. Suppose that 
\begin{align}\label{eqq17}
\bigg|\frac{1}{\log \omega(x)}\sum_{\frac{x}{\omega(x)}\leq n\leq x}\frac{g_1(n)g_2(n+h)}{n}\bigg|>\varepsilon^2.    
\end{align}
Let $\exp_2(\varepsilon^{-1})\leq y\leq \log \log x$ be arbitrary. Then there exists a set $\mathcal{N}\subset [1,y]$ such that for all $m\in \mathcal{N}$ we have
\begin{align*}
\sum_{\substack{2^m\leq p<2^{m+1}\\|g_1(p)g_2(p)|>\frac{1}{2}}}1\geq (1-\varepsilon)\cdot \frac{2^m}{m\log 2},    
\end{align*}
with $\mathcal{N}$ being large in the sense that
\begin{align*}
\frac{1}{\log y}\sum_{\substack{n\leq y\\n\in \mathcal{N}}}\frac{1}{n}\geq 1-\varepsilon.    
\end{align*}
\end{lemma}

\textbf{Proof.} Suppose for the sake of contradiction that such a set $\mathcal{N}$ does not exist. Then by the prime number theorem we have
\begin{align}\label{eqq44}
\sum_{\substack{2^m\leq p<2^{m+1}\\|g_1(p)g_2(p)|\leq \frac{1}{2}}}1\geq \frac{\varepsilon}{2}\cdot \frac{2^m}{m\log 2}    
\end{align}
for all $m\in \mathcal{N}_1\subset [1,y]$ with $\mathcal{N}_1$ being a set with the property
\begin{align}\label{eqq45}
\sum_{m\in \mathcal{N}_1}\frac{1}{m}\geq \frac{\varepsilon}{2}\log y.
\end{align}
In particular, from \eqref{eqq44} we have
\begin{align*}
 \sum_{\substack{2^m\leq p<2^{m+1}\\|g_1(p)g_2(p)|\leq \frac{1}{2}}}\frac{1}{p}\geq \frac{\varepsilon}{8m}  
\end{align*}
for $m\in \mathcal{N}_1$. Summing over $m\in \mathcal{N}_1$ and using \eqref{eqq45}, we conclude that
\begin{align*}
 \sum_{\substack{p\leq 2^{y+1}\\|g_1(p)g_2(p)|\leq \frac{1}{2}}}\frac{1}{p}\geq \frac{\varepsilon^2}{16}\log y.  
\end{align*}
Hence, for at least one of $j=1$ and $j=2$ we have
\begin{align}\label{eqq46}
\sum_{\substack{p\leq 2^{y+1}\\|g_j(p)|\leq \frac{1}{\sqrt{2}}}}\frac{1}{p}\geq \frac{\varepsilon^2}{32}\log y.    
\end{align}
Fix such $j\in \{1,2\}$. Let 
\begin{align*}
\mathcal{P}:=\{\varepsilon^{-10}\leq p\leq 2^{y+1}:|g_j(p)|\leq \frac{1}{\sqrt{2}}\},    
\end{align*}
and let $\mu^2_{\mathcal{P}}(n)$ be the indicator function of integers $n$ that are not divisible by $p^2$ for any $p\in \mathcal{P}$. Note that if $\mu^2_{\mathcal{P}}(n)=1$, then 
\begin{align*}
|g_j(n)|\leq \left(\frac{1}{\sqrt{2}}\right)^{\omega_{\mathcal{P}}(n)},   
\end{align*}
where $\omega_{\mathcal{P}}(n)$ is the number of prime factors of $n$ from $\mathcal{P}$. In particular, we have $|g_j(n)|\leq \varepsilon^{10}$ whenever $\omega_{\mathcal{P}}(n)\geq \varepsilon^{-1}$ (and still $\mu^2_{\mathcal{P}}(n)=1$). In conclusion, if we show that
\begin{align}\label{eqq18}
\frac{1}{\log \omega(x)}\sum_{\substack{\frac{x}{\omega(x)}\leq n\leq x\\\mu^2_{\mathcal{P}}(n)=0\\\textnormal{or}\,\, \omega_{\mathcal{P}}(n)<\varepsilon^{-1}}}\frac{1}{n}\leq \varepsilon^3,    
\end{align}
then \eqref{eqq17} is violated, giving the desired contradiction. We are now left with showing \eqref{eqq18}, and for this we use some basic sieve theory. Note that
\begin{align*}
\frac{1}{\log \omega(x)}\sum_{\substack{\frac{x}{\omega(x)}\leq n\leq x\\\mu^2_{\mathcal{P}}(n)=0}} \frac{1}{n}\leq \sum_{p\in \mathcal{P}}\frac{1}{\log \omega(x)}\sum_{\frac{x}{p^2\omega(x)}\leq m\leq \frac{x}{p^2}}\frac{1}{p^2m}\ll\sum_{p\in \mathcal{P}}\frac{1}{p^2}\ll \varepsilon^{10}.    
\end{align*}
Note also that if $\omega_{\mathcal{P}}(n)=M$ and $\mu_{\mathcal{P}}^2(n)=1$, then we may write $n=p_1\cdots p_Mm$ with $p_i\in \mathcal{P}$ and $\omega_{\mathcal{P}}(m)=0$. Hence, by the sieve of Eratosthenes and Mertens' theorem,
\begin{align*}
&\sum_{\substack{\frac{x}{\omega(x)}\leq n\leq x\\\omega_{\mathcal{P}}(n)<\varepsilon^{-1}\\\mu_{\mathcal{P}}^2(n)=1}}\frac{1}{n}\leq \varepsilon^{-1}\max_{M<\varepsilon^{-1}}\sum_{p_1,\ldots, p_M\leq 2^{y+1}}\sum_{\substack{\frac{x}{\omega(x)p_1\cdots p_M}\leq m\leq \frac{x}{p_1\cdots p_M}\\\omega_{\mathcal{P}}(m)=0}}\frac{1}{p_1\cdots p_M m}\\
&\ll \varepsilon^{-1}\max_{M<\varepsilon^{-1}}\sum_{p_1,\ldots, p_M\leq 2^{y+1}}\frac{1}{p_1\cdots p_M}\prod_{p\in \mathcal{P}}\left(1-\frac{1}{p}\right)\cdot \log \omega(x)\\
&\ll \varepsilon^{-1}\max_{M<\varepsilon^{-1}}\sum_{p_1,\ldots, p_M\leq 2^{y+1}}\frac{1}{p_1\cdots p_M}\exp\left(-\sum_{p\in \mathcal{P}}\frac{1}{p}\right)\cdot \log \omega(x)\\
&\ll \varepsilon^{-1}(\log y)^{\varepsilon^{-1}}y^{-\frac{\varepsilon^2}{32}}\log \omega(x)\ll \varepsilon^{10}\log \omega(x)
\end{align*}
by \eqref{eqq46} and the fact that $y\geq \exp_2(\varepsilon^{-1})$. Combining the above estimates, we obtain \eqref{eqq18}, and hence also the statement of the lemma.\qedd\\

\begin{lemma}[Dealing with pretentious functions]\label{le_C} Let the notations be as in Theorem \ref{theo_bincorr}. Suppose that $|\delta_1|>\varepsilon^2$ and $|\delta_2|>\varepsilon^2$, where $\delta_1$ and $\delta_2$ are as in \eqref{eqq82}. Let $\exp_2(\varepsilon^{-1})\leq y\leq \log \log x$ be arbitrary. Then there exists a set $\mathcal{N}'\subset [1,y]$ such that for all $m\in \mathcal{N}'$ we have
\begin{align*}
\frac{m\log 2}{2^m}\bigg|\sum_{2^m\leq p<2^{m+1}}(1-g_1(p)g_2(p))\bigg|<\varepsilon,
\end{align*}
with $\mathcal{N}'$ large in the sense that
\begin{align*}
\frac{1}{\log y}\sum_{\substack{n\leq y\\n\in \mathcal{N}'}}\frac{1}{n}\geq 1-\varepsilon.
\end{align*}
\end{lemma}

\textbf{Proof.} Note that $1-g_1(p)g_2(p)\geq 0$ always holds. Arguing just as in the proof of Lemma \ref{le_A}, we see that if the statement failed, we would have
\begin{align*}
 \sum_{p\leq 2^{y+1}}\frac{1-g_1(p)g_2(p)}{p}\geq \frac{\varepsilon^2}{8}\log y.  
\end{align*}
In particular, by the inequality $(1-a)+(1-b)\geq 1-ab$ for $a,b\in [-1,1]$, for at least one of $j=1$ and $j=2$ we would have
\begin{align}\label{eqq85}
\sum_{p\leq 2^{y+1}}\frac{1-g_j(p)}{p}\geq \frac{\varepsilon^2}{16}\log y.    
\end{align}
Now, by \eqref{eqq85} and a version of Halász's theorem for real-valued multiplicative functions \cite{hall1}, we have
\begin{align*}
|\delta_j|=\left|\frac{1}{x}\sum_{n\leq x}g_j(n)\right|\ll \exp\left(-\frac{1}{10}\sum_{p\leq 2^{y+1}}\frac{1-g_j(p)}{p}\right)\ll \exp\left(-\frac{\varepsilon^2}{200}\log y\right)\ll \varepsilon^{10} \end{align*}
for $y\geq \exp_2(\varepsilon^{-1})$, and this contradicts $|\delta_j|>\varepsilon^2$, proving the lemma.\qedd\\

Now we return to \eqref{eqq5} and consider two cases separately. Suppose first that $|\delta_1|, |\delta_2|>\varepsilon^2$. Let $y=\exp_2(\varepsilon^{-1})^2$. Then, if $\mathcal{N}'\subset [1,y]$ is the set in Lemma \ref{le_C} and $\mathcal{M}$ is the set in Lemma \ref{le_entropy} (which is independent of $c_p$), taking $c_p=1$ we deduce from \eqref{eqq5} and Lemma \ref{le_C} that
\begin{align}\label{eqq19}\begin{split}
f_{x,\omega}(h)&=\frac{m\log 2}{2^m}\sum_{2^m\leq p<2^{m+1}}g_1(p)g_2(p)f_{x,\omega}(h)+o_{\varepsilon\to 0}(1)\\
&=\delta_1\delta_2+\frac{m\log 2}{2^m}\sum_{2^m\leq p<2^{m+1}}\tilde{f}_{x,\omega}(ph)+o_{\varepsilon\to 0}(1)
\end{split}
\end{align}
for $m\in \mathcal{M}\cap \mathcal{N}'$. We can pick some $m\in \mathcal{M}\cap \mathcal{N}'$ with $m\in [\sqrt{y},y]$, since we have the lower bound
\begin{align*}
\sum_{\substack{m\in \mathcal{M}\cap \mathcal{N}'\\m\in [\sqrt{y},y]}}\frac{1}{m}\geq \log y-\frac{1}{2}\log y-\varepsilon^{-100}- \varepsilon \log y\geq \frac{1}{3}\log y \end{align*}
for $y=\exp_2(\varepsilon^{-1})^2$.\\

Consider then  the case where either $|\delta_1|\leq \varepsilon^{2}$ or $|\delta_2|\leq \varepsilon^{2}$. We may suppose that \eqref{eqq17} holds, since otherwise Theorem \ref{theo_bincorr} holds by the fact that $\delta_1\delta_2+o_{\varepsilon\to 0}(1)=o_{\varepsilon\to 0}(1)$ in this situation.  Let $y=\exp_2(\varepsilon^{-1})^2$. Taking $m\in \mathcal{M}\cap \mathcal{N}$ (with $\mathcal{N}\subset [1,y]$ as in Lemma \ref{le_A}) and $c_p=g_1(p)^{-1}g_2(p)^{-1}1_{|g_1(p)g_2(p)|\geq \frac{1}{2}}$ in \eqref{eqq5}, we see from Lemma \ref{le_A} that
\begin{align*}
&f_{x,\omega}(h)+o_{\varepsilon\to 0}(1)\\
&=\delta_1\delta_2\frac{m\log 2}{2^m} \sum_{\substack{2^m\leq p<2^{m+1}\\|g_1(p)g_2(p)|\geq \frac{1}{2}}}(g_1(p)g_2(p))^{-1}+\frac{m\log 2}{2^m}\sum_{\substack{2^m\leq p<2^{m+1}\\|g_1(p)g_2(p)|\geq \frac{1}{2}}}(g_1(p)g_2(p))^{-1}\tilde{f}_{x,\omega}(ph).    \end{align*}
for $m\in \mathcal{M}\cap \mathcal{N}$, which again contains an element $m\in [\sqrt{y},y]$ by the same argument as above. We know that $|(g_1(p)g_2(p))^{-1}|\leq 2$ for all $2^m\leq p<2^{m+1}$, except for at most $10\varepsilon\frac{2^m}{m}$ exceptions. Since by assumption $\delta_1\delta_2=O(\varepsilon)$, we deduce that 
\begin{align*}
f_{x,\omega}(h)&=\delta_1\delta_2+\frac{m\log 2}{2^m}\sum_{2^m\leq p<2^{m+1}}2a_p\tilde{f}_{x,\omega}(ph)+o_{\varepsilon\to 0}(1)
\end{align*}
for $m\in \mathcal{M}\cap \mathcal{N}$, where $a_p:=\frac{1}{2}(g_1(p)g_2(p))^{-1}1_{|g_1(p)g_2(p))|\geq \frac{1}{2}}$. In conclusion, regardless of the values of $\delta_j$, Theorem \ref{theo_bincorr} will follow once we prove that
\begin{align}\label{eqq20}
\frac{m}{2^m \log \omega(x)}\sum_{2^m\leq p<2^{m+1}}a_p \sum_{\frac{x}{\omega(x)}\leq n\leq x}\frac{(g_1(n)-\delta_1)(g_2(n+ph)-\delta_2)}{n}= o_{\varepsilon\to 0}(1)
\end{align}
for arbitrary $a_p\in \mathbb{D}$ and $m\in [\exp_2(\varepsilon^{-1}),\exp_2(\varepsilon^{-1})^2]$.

\section{Circle method estimates}\label{sec:circle}

We proceed to prove \eqref{eqq20} by applying the circle method and (slightly modified versions of) the short exponential sum estimates for multiplicative functions due to Matom\"aki-Radziwi\l{}\l{} and Tao \cite{mrt}. We start with two lemmas, the first of which reduces \eqref{eqq20} to bounding a short exponential sum and the second of which shows that the set of large frequencies of the exponential sum has small cardinality.

\begin{lemma}[A circle method estimate] \label{le_circle}
Let $\eta>\varepsilon>0$ be small, $h$ an integer with $1\leq |h|\leq \varepsilon^{-1}$, and $\exp_2(\varepsilon^{-1})\leq H\leq \log y$. For any complex numbers $a_p\in \mathbb{D}$, introduce the exponential sum
\begin{align*}
S_H(\theta):=\sum_{P\leq p<2P}a_p e(p\theta),    
\end{align*}
where $P:=\varepsilon^{10}H$. Let $\Xi_H$ be the set of residue classes $\xi \in \mathbb{Z}/H\mathbb{Z}$ that satisfy
\begin{align}\label{eqq91}
\left|S_H\left(-\frac{h\xi}{H}\right)\right|\geq \eta^{2}\frac{P}{\log H}.    
\end{align}
Then, for any functions $g_1',g_2':\mathbb{N}\to \mathbb{C}$ with $|g_1'(n)|, |g_2'(n)|\leq 2$, we have
\begin{align*}
 &\bigg|\frac{\log P}{P}\sum_{P\leq p<2P}a_p\sum_{y\leq n\leq y+H}g_1'(n)g_2'(n+ph)\bigg|\leq \eta H+10\sum_{\xi \in \Xi_H}\bigg|\sum_{y\leq n\leq y+H}g_1'(n)e\left(-\frac{\xi n}{H}\right)\bigg|.
\end{align*}
\end{lemma}

\textbf{Proof.} This follows from \cite[Lemma 3.6]{tao}, writing it using different notation.\qedd\\

In order to make use of Lemma \ref{le_circle}, we must know that the exceptional set $\Xi_H$ in that theorem is not too large. Indeed, we have the following bound.

\begin{lemma}[Cardinality of large Fourier coefficients] \label{le_exceptional} Let the notations be as in Lemma \ref{le_circle}, and assume that $H$ is a prime. Then we have $|\Xi_H|\ll \eta^{-20}$.
\end{lemma}

\textbf{Proof.} Since $1\leq |h|\leq \varepsilon^{-1}$ and $H$ is a prime, the number of those $\xi$ that satisfy \eqref{eqq91} remains unchanged when $h$ is replaced by $1$ in that formula. In \cite[Lemma 3.7]{tao}, it was proved using a fourth moment bound and the Selberg sieve that $|S_{H}(-\frac{\xi}{H})|\geq \frac{\eta^2 P}{\log H}$for $\ll_{\eta} 1$ values of $\xi\in \mathbb{Z}/H\mathbb{Z}$, but the same proof gives the claimed quantitative bound.\qedd\\

To make use of the two lemmas above, we split in \eqref{eqq20} the sum over $n$ into sums of length $H$, where $H$ is a prime belonging to $[\varepsilon^{-10}\cdot 2^m,2\varepsilon^{-10}\cdot 2^m]$, and approximate the sum with an integral, after which \eqref{eqq20} is reduced to
\begin{align}\label{eqq38}
\frac{1}{\log \omega(x)}\int_{\frac{x}{\omega(x)}}^x \frac{m}{2^m }\sum_{2^m\leq p<2^{m+1}} \frac{a_p}{H}\sum_{y\leq n\leq y+H}(g_1(n)-\delta_1)(g_2(n+ph)-\delta_2)\, \frac{dy}{y}=o_{\varepsilon\to 0}(1).    
\end{align}
By Lemmas \ref{le_circle} and \ref{le_exceptional}, it suffices to show that
\begin{align}\label{eqq39}
\sup_{\alpha \in \mathbb{R}}\frac{1}{\log \omega(x)}\int_{\frac{x}{\omega(x)}}^x \frac{1}{y}\bigg|\frac{1}{H}\sum_{y\leq n\leq y+H}(g_1(n)-\delta_1)e(\alpha n)\bigg| \, dy=o_{\varepsilon\to 0}(1),    \end{align}
for $H\in [\exp_3(\frac{1}{2}\varepsilon^{-1}),\exp_3(2\varepsilon^{-1})]$, where $\exp_3$ is the third iterated exponential. Indeed, if the left-hand side of \eqref{eqq39} is $\leq F(\varepsilon)$, where $F(u)\to 0$ as $u\to 0$ is a slowly decaying function, one can take $\eta=F(\varepsilon)^{0.01}$ in Lemma \ref{le_exceptional} to deduce \eqref{eqq38}. Covering the interval $[\frac{x}{\omega(x)},x]$ with dyadic intervals, \eqref{eqq39} will follow from 
\begin{align}\label{eqq41}
\sup_{\alpha \in \mathbb{R}}\frac{1}{X}\int_{X}^{2X} \bigg|\frac{1}{H}\sum_{y\leq n\leq y+H}(g_1(n)-\delta_1)e(\alpha n)\bigg| \, dy=o_{\varepsilon\to 0}(1)
\end{align}
for all $X\in [\frac{x}{\omega(x)},\frac{x}{2}]$ and all $H\in [\exp_3(\frac{1}{2}\varepsilon^{-1}),\exp_3(2\varepsilon^{-1})]$. This is what we set out to prove, following \cite{mrt}.\\

It is natural to split the supremum over $\alpha$ in \eqref{eqq41} to major and minor arcs, defined using Dirichlet's approximation theorem as
\begin{align}\label{eqtn25}\begin{split}
\mathfrak{M}&:=\left\{\theta \in \mathbb{R}:\,\, \left|\theta-\frac{a}{q}\right|\leq \frac{W}{qH}\,\,\textnormal{with}\,\, a\in \mathbb{Z},\,\, q< W,\,\, (a,q)=1\right\}\quad \textnormal{and}\\
\mathfrak{m}&:=\mathbb{R}\setminus{\mathfrak{M}}\subset \left\{\theta \in \mathbb{R}:\,\, \left|\theta-\frac{a}{q}\right|\leq \frac{W}{qH}\,\, \textnormal{with}\,\, a\in \mathbb{Z},\,\, q\in \left[W,\frac{H}{W}\right],\,\, (a,q)=1\right\},\,\, \textnormal{with}\\
W&:=\log^5 H\leq \exp(5\exp(2\varepsilon^{-1})).
\end{split}
\end{align}

In the case of the major arcs, the exponential $e(\alpha n)$ can essentially be replaced with $e(\frac{an}{q})$, and this will lead us to study the distribution of the multiplicative function $g_1$ in arithmetic progressions over short intervals. For that purpose, we prove a lemma that is closely related to \cite[Theorem 1]{matomaki-radziwill} and \cite[Theorem A.1]{mrt}. For this lemma, we need to introduce the same ``nicely factorable'' set as in \cite[Section 2]{matomaki-radziwill} and \cite[Definition 2.1]{mrt}. 

\begin{definition}\label{def3} Let $10<P_1<Q_1\leq X$ and $\sqrt{X}\leq X_0\leq X$, with $Q_1\leq \exp(\sqrt{\log X_0})$. For $j>1$, set
\begin{align*}
P_j:=\exp(j^{4j}(\log Q_1)^{j-1}\log P_1),\quad Q_j:=\exp(j^{4j+2}(\log Q_1)^{j}).    
\end{align*}
Letting $J$ be the largest integer such that $Q_J\leq \exp(\sqrt{\log X_0})$, we define $\mathcal{S}_{P_1,Q_1,X_0,X}$ as the set of those $1\leq n\leq X$ that have at least one prime factor from each of the intervals $[P_j,Q_j]$ for all $1\leq j\leq J$.
\end{definition}

For a specific choice of the parameters, present in the next lemma, we denote
\begin{align}\label{eqq50}
\mathcal{S}:=\mathcal{S}_{P_1,Q_1,X_0,X},\quad \textnormal{where}\quad P_1=W^{200},\,\, Q_1=\frac{H}{W^3},\,\, X_0=\sqrt{X}.     
\end{align}

\begin{lemma}[Uniform distribution of multiplicative functions in short intervals]\label{le_mr} Let $\varepsilon>0$ be small, $X\geq 100$ large, and $H\in [\exp_2(\frac{1}{10}\varepsilon^{-1}), \log \log X]$. Let $g:\mathbb{N}\to [-1,1]$ be a real-valued multiplicative function.  Further, let  $b,q\in \mathbb{N}$ with $1\leq b\leq q\leq W\in [\log^5 H,\log^{10} H]$. Then, if $\mathcal{S}$ is as in \eqref{eqq50}, we have
\begin{align}\label{eqq86}
\frac{1}{X}\int_{X}^{2X}\bigg|\frac{1}{H}\sum_{\substack{y\leq n\leq y+H\\n\equiv b\pmod q}}g(n)1_{\mathcal{S}}(n)-\frac{1}{X}\sum_{\substack{X\leq n\leq 2X\\n\equiv b\pmod q}}g(n)1_{\mathcal{S}}(n)\bigg|\,dy\ll W^{-10}.
\end{align}
\end{lemma}

\begin{remark} If the bound on the right-hand side of \eqref{eqq86} was replaced with $W^{-0.001}$, the proof of the lemma would work even when $g(n)1_{\mathcal{S}}(n)$ is replaced with $g(n)$. However, for larger values of $q$ we need to introduce the nicely factorable set $\mathcal{S}$ to get better error terms.
\end{remark}

\textbf{Proof of Lemma \ref{le_mr}.} We will first reduce to primitive residue classes $b\pmod q$. Let $d_0=(b,q)$, $b_0=\frac{b}{(b,q)}$ and $q_0=\frac{q}{(b,q)}$. Then we have
\begin{align}\label{eqq24}
\frac{1}{H}\sum_{\substack{y\leq n\leq y+H\\n\equiv b\pmod q}}g(n)1_{\mathcal{S}}(n)=\frac{1}{H}\sum_{\substack{y/d_0\leq n'\leq (y+H)/d_0\\n\equiv b_0\pmod{q_0}}}g(d_0n')1_{\mathcal{S}}(n'),
\end{align}
since $1_{\mathcal{S}}(d_0n')=1_{\mathcal{S}}(n')$ for $d_0\leq q\leq W <P_1$. Since the residue class $b_0\pmod{q_0}$ is primitive, we may use a Dirichlet character expansion to write the right-hand side of \eqref{eqq24} as
\begin{align}\label{eqq31}\begin{split}
&\frac{1}{H}\frac{1}{\varphi(q_0)}\sum_{\chi \pmod{q_0}}\bar{\chi}(b_0)\sum_{y/d_0\leq n'\leq (y+H)/d_0}g(d_0n')1_{\mathcal{S}}(n')\chi(n')\\
&=\frac{1}{H}\frac{1}{\varphi(q_0)}\sum_{\chi \pmod{q_0}}\bar{\chi}(b_0)\sum_{t\mid d_0^{\infty}}\,\sum_{\substack{y/d_0\leq n'\leq (y+H)/d_0\\t\mid n'\\(\frac{n'}{t},d_0)=1}}g(d_0n')1_{\mathcal{S}}(n')\chi(n'),
\end{split}
\end{align} 
where $t\mid d_0^{\infty}$ means that $t\mid d_0^k$ for some $k$. Since we have the condition $(\frac{n'}{t},d_0)=1$, we may use multiplicativity to write this as
\begin{align}\label{eqq25}
\frac{1}{\varphi(q_0)}\sum_{\chi \pmod{q_0}}\bar{\chi}(b_0)\sum_{t\mid d_0^{\infty}}\frac{g(d_0t)\chi(t)}{d_0t}\frac{d_0t}{H}\sum_{y/(d_0t)\leq m\leq (y+H)/(d_0t)}g(m)1_{\mathcal{S}}(m)\chi\psi_0(m),
\end{align}
where $\psi_0(m)=1_{(m,d_0)=1}$ is the principal character $\pmod{d_0}$ and we used the fact that $1_{\mathcal{S}}(tm)=1_{\mathcal{S}}(m)$ for $t$ having no prime factors that are larger than $d_0\leq q\leq W<P_1$. By crude estimation, the contribution of the terms $t\geq H^{\varepsilon}$ to \eqref{eqq25} is $\ll H^{-\varepsilon}$, so we may assume that $t<H^{\varepsilon}$. We now wish to compare the short sums in \eqref{eqq25} to the corresponding long sums.\\

Suppose first that $\chi$ is real-valued. Then we may apply the Matom\"aki-Radziwi\l{}\l{} theorem \cite[Theorem 3]{matomaki-radziwill} to the real-valued multiplicative function $g\chi\psi_0$ conclude that
\begin{align}\label{eqq26}\begin{split}
&\frac{d_0t}{H}\sum_{y/(d_0t)\leq m\leq (y+H)/(d_0t)}g(m)1_{\mathcal{S}}(m)\chi\psi_0(m)\\
&=\frac{d_0t}{X}\sum_{X/(d_0t)\leq m\leq 2X/(d_0t)}g(m)1_{\mathcal{S}}(m)\chi\psi_0(m)+E_{\chi,H}(y),
\end{split}
\end{align} 
for $y\in [X,2X]$, with the error $E_{\chi, H}(y)$ satisfying the $L^2$ bound
\begin{align}\label{eqq27}
\frac{1}{X}\int_{X}^{2X}|E_{\chi, H}(y)|^2\,dy\ll \frac{(\log H)^{\frac{1}{3}}}{P_1^{\frac{1}{10}}}+(\log X)^{-\frac{1}{50}}\ll W^{-19},
\end{align}
since $W\in [\log^5 H, \log^{10} H]$, and $P_1=W^{200}$ in our definition of $\mathcal{S}$.\\

Suppose then that $\chi$ is complex-valued. We again write \eqref{eqq26}, and want to obtain an $L^2$ bound for the error $E_{\chi,H}(y)$. By an argument of Granville and Soundararajan (see \cite[Lemma C.1]{mrt}), the fact that $g\psi_0$ is real and $\chi$ is complex (and that $q\leq (\log_3 X)^{10}$) leads to
\begin{align}\label{eqq28}
\inf_{|t|\leq x}\mathbb{D}(g\chi\psi_0,n^{it};x)\geq \frac{1}{10}\sqrt{\log \log x}.
\end{align}
Now we appeal to a variant of the Matom\"aki-Radziwi\l{}\l{} theorem, established by Matom\"aki, Radziwi\l{}\l{} and Tao in \cite[Theorem A.2]{mrt}. This result (applied with $h=H$ and $h=X$ separately) gives
\begin{align}\label{eqq51}
\frac{1}{X}\int_{X}^{2X}|E_{\chi, H}(y)|^2\,dy\ll \exp\left(-\inf_{|t|\leq X}\frac{\mathbb{D}(g\chi\psi_0,n^{it};X)^2}{2}\right)+\frac{(\log H)^{\frac{1}{3}}}{P_1^{\frac{1}{10}}}+(\log X)^{-\frac{1}{50}}\ll W^{-19}  
\end{align}
by \eqref{eqq28}.\\

Now, for all characters $\chi \pmod{q_0}$, we have \eqref{eqq26} with the error bound \eqref{eqq27}. Note also that $\sum_{t\mid d_0^{\infty}}\frac{1}{d_0t}=\frac{1}{d_0}\prod_{p\mid d_0}(1+\frac{1}{p}+\frac{1}{p^2}+\cdots)\ll \frac{\log d_0}{d_0}$. Hence, applying the triangle inequality, and summing over $\chi$ and $t\mid d_0^{\infty}$, we see that \eqref{eqq25} equals
\begin{align}\label{eqq30}
\frac{1}{\varphi(q_0)}\sum_{\chi \pmod{q_0}}\bar{\chi}(b_0)\sum_{t\mid d_0^{\infty}}\frac{g(d_0t)\chi(t)}{d_0t}\frac{d_0t}{X}\sum_{X/(d_0t)\leq m\leq 2X/(d_0t)}g(m)1_{\mathcal{S}}(m)\chi\psi_0(m)+E(y),
\end{align}
with the error term $E(y)$ satisfying
\begin{align*}
\frac{1}{X}\int_{X}^{2X}|E(y)|^2\, dy\ll W^{-10}.
\end{align*}
We can then reverse the deduction that led to \eqref{eqq25} to conclude that \eqref{eqq30} (and hence \eqref{eqq24}) equals
\begin{align*}
\frac{1}{X}\sum_{\substack{X\leq n\leq 2X\\n\equiv b\pmod q}}g(n)1_{\mathcal{S}}(n)+E(y).
\end{align*}
This completes the proof.\qedd\\

The major arc case $\alpha \in \mathfrak{M}$ of \eqref{eqq41} is dealt with the following Lemma, whose proof uses Lemma \ref{le_mr} as an ingredient.

\begin{lemma}[Major arc estimate]\label{le_major} Let $\varepsilon>0$ be small, $x\geq 100$ large, $\omega(X)$ as in Theorem \ref{theo_bincorr}, and $H\in [\exp_3(\frac{1}{2}\varepsilon^{-1}), \exp_3(2\varepsilon^{-1})]$. Let $g_1:\mathbb{N}\to [-1,1]$ be a multiplicative function satisfying $g_1\in \mathcal{U}_{\omega}(x,\exp_2(\varepsilon^{-2}),2/\exp_2(\varepsilon^{-2}),\delta_1)$.  Then we have
\begin{align*}
\sup_{\alpha \in \mathfrak{M}}\frac{1}{X}\int_{X}^{2X}\left|\frac{1}{H}\sum_{y\leq n\leq y+H}(g_1(n)-\delta_1)e(\alpha n)\right|\,dy=o_{\varepsilon\to 0}(1)
\end{align*}
for all $X\in [\frac{x}{\omega(x)},\frac{x}{2}]$, with the major arcs $\mathfrak{M}$ as in \eqref{eqtn25}.
\end{lemma}

\textbf{Proof.} This is proved in Appendix \ref{a:b}.\qedd\\

The minor arc case $\alpha \in \mathfrak{m}$ of \eqref{eqq41}, in turn, is taken care of by the next lemma.

\begin{lemma}[Minor arc estimate] \label{le_minor} Let $\varepsilon>0$ be small, $x\geq 100$ large, and suppose that $H\in [\exp_3(\frac{1}{2}\varepsilon^{-1}),\log \log x]$. Then, for any multiplicative function $g_1:\mathbb{N}\to [-1,1]$ and for any $\delta_1\in [-1,1]$ we have
\begin{align}\label{eqq33}
\sup_{\alpha \in \mathfrak{m}}\frac{1}{X}\int_{X}^{2X}\left|\frac{1}{H}\sum_{y\leq n\leq y+H}(g_1(n)-\delta_1)e(\alpha n)\right|\, dy\ll (\log H)^{-\frac{1}{10}}
\end{align}
for all $X\in [\sqrt{x},x]$, with the minor arcs $\mathfrak{m}$ as in \eqref{eqtn25}.
\end{lemma}

\textbf{Proof.} This is proved in Appendix \ref{a:b}.\qedd\\

With these lemmas available, Theorem \ref{theo_bincorr} quickly follows.\\

\textbf{Proof of Theorem \ref{theo_bincorr}.} We reduced the proof of the theorem to proving \eqref{eqq41}. As was observed after Lemma \ref{le_stability}, we may assume that $g_1\in \mathcal{U}_{\omega}(x,\exp_2(\varepsilon^{-2}),2/\exp_2(-\varepsilon^{-2}),\delta_1)$. Now, if $\alpha\in \mathfrak{M}$ in the supremum present in that formula, we appeal to Lemma \ref{le_major}. In the opposite case $\alpha \in \mathfrak{m}$, we appeal to Lemma \ref{le_minor}. In both cases, we get a bound of $o_{\varepsilon\to 0}(1)$ for the left-hand side of \eqref{eqq41}. This finishes the proof.\qedd

\section{Proofs of the applications}\label{sec:apps}

\textbf{Proof of Theorem \ref{theo_omega}.} Given any real numbers $z,w\in [-1,1]$, define the multiplicative functions $g_1,g_2:\mathbb{N}\to [-1,1]$ by setting at prime powers
\begin{align*}
g_1(p^j)=\begin{cases}1\,\, \textnormal{if}\,\, p\leq x^{a}\\z\,\, \textnormal{if}\,\,p>x^{a},\end{cases}\quad g_2(p^j)=\begin{cases}1\,\,\textnormal{if}\,\, p\leq x^{b}\\w\,\,\textnormal{if}\,\, p>x^{b}.\end{cases}
\end{align*}
 We will apply Theorem \ref{theo_bincorr} to $g_1$ and $g_2$, and then use a generating function argument to deduce Theorem \ref{theo_omega}. In order to use Theorem \ref{theo_bincorr}, we must verify that $g_1\in \mathcal{U}(x,\varepsilon^{-1},\varepsilon)$ for all $x\geq x_0(\varepsilon)$.\\

First observe that $g_1(n)=z^{\omega_{>x^{a}}(n)}$, so for any $c,q\in \mathbb{N}$ we have
\begin{align*}
\frac{1}{x}\sum_{\substack{x\leq n\leq 2x\\n\equiv c\pmod{q}}}g_1(n)&=\sum_{0\leq k< \frac{1}{a}}z^k\cdot \frac{1}{x}\sum_{\substack{x\leq n\leq 2x\\n\equiv c\pmod{q}}}1_{\omega_{>x^{a}}(n)=k}.
\end{align*}
From this we see that $g_1\in \mathcal{U}(x,\varepsilon^{-1},\varepsilon)$ for all $x\geq x_0(\varepsilon)$ will follow, once we show that
\begin{align*}
\frac{1}{x}\sum_{\substack{x\leq n\leq 2x\\n\equiv c\pmod{q}}}1_{\omega_{>x^{a}}(n)=k}=\frac{1}{qx}\sum_{x\leq n\leq 2x}1_{\omega_{>x^{a}}(n)=k}+o_q(1)   
\end{align*}
as $x\to \infty$ for all fixed $c,q,k\in \mathbb{N}$. Write $d_0=(c,q)$, $c'=\frac{c}{d_0}$, $q'=\frac{q}{d_0}$. Then we have
\begin{align}\label{eqq66}
\frac{1}{x}\sum_{\substack{x\leq n\leq 2x\\n\equiv c\pmod{q}}}1_{\omega_{>x^{a}}(n)=k}&=\frac{1}{x}\sum_{\substack{\frac{x}{d_0}\leq n'\leq \frac{2x}{d_0}\\n'\equiv c'\pmod{q'}}}1_{\omega_{>x^{a}}(n')=k}:=S_k,
\end{align}
because $\omega_{>x^{a}}(d_0n')=\omega_{>x^{a}}(n')$ for all $d_0<x^{a}$. Let $b^{-1}\pmod{q}$ denote the inverse of $b$ modulo $q$. Using the fact that $1_{\omega_{>x^{a}}(n)=0}=1_{P^{+}(n)\leq x^{a}}$, we have
\begin{align*}
S_k=\sum_{\substack{x^{a}<p_1<\ldots<p_k\leq x\\p_1\cdots p_k\leq x}}\frac{1}{x}\sum_{\substack{\frac{x}{d_0p_1\cdots p_k}\leq m\leq \frac{2x}{d_0p_1\cdots p_k}\\m\equiv c'(p_1\cdots p_k)^{-1}\pmod{q'}}}1_{P^{+}(m)\leq x^{a}}+o_{q'}(1),
\end{align*}
with the $o(1)$ term coming from those numbers $n'\leq x$ such that $p^2\mid n'$ for some $p>x^{a}$. As is well-known, smooth numbers are uniformly distributed in arithmetic progressions to fixed moduli (see for instance \cite[Formula (6.1)]{ht}), in the sense that
\begin{align}\label{eqq90}
\frac{1}{y}|\{y\leq n\leq 2y:\,\, P^{+}(n)\leq y^{u}, n\equiv c \pmod{q'}\}|=\frac{1}{q'}\rho\left(\frac{1}{u}\right)+o_{q'}(1),  
\end{align}
for $u\in [0,1]$ and $y\to \infty$, with $\rho(\cdot)$ being the Dickmann function. Therefore, 
\begin{align}\label{eqq87}
S_k=\frac{1}{q'd_0}\sum_{\substack{x^{a}<p_1<\ldots<p_k\leq x\\p_1\cdots p_k\leq x}}\frac{1}{p_1\cdots p_k}\left(\rho\left(\frac{\log \frac{x}{d_0p_1\cdots p_k}}{a\log x}\right)+o_{q'}(1)\right).
\end{align}
One easily sees that $x\mapsto \rho(x)$ is a Lipschitz function, so that $|\rho(u)-\rho(v)|\leq C|u-v|$ for all $u,v\geq 0$ with some constant $C>0$. Hence, we can use the prime number theorem in the form that the $n$th prime is asymptotic to $n\log n$ and approximate the term involving $\rho(\cdot)$ in \eqref{eqq87} to deduce that
\begin{align}\label{eqq88}
S_k=\frac{1}{q'd_0k!}\sum_{\substack{x^{a}<n_1,\ldots,n_k\leq x\\n_1\cdots n_k\leq x}}\frac{1}{n_1\cdots n_k(\log n_1)\cdots (\log n_k)}\rho\left(\frac{\log \frac{x}{n_1\cdots n_k}}{a\log x}\right)+o_{q'}(1).    
\end{align}
Here we have estimated trivially as $o_{q'}(1)$ the contribution of the tuples $(n_1,\ldots, n_k)$ with two of the $n_i$ equal, or with $n_i\in [\frac{x^{a}}{2\log x},x^{a}]\cup [\frac{x}{2\log x},x]$ for some $i$, as for them it is not necessarily the case that the $n_i$th prime belongs to $[x^{a},x]$. Approximating the expression \eqref{eqq88} with an integral, again using the fact that $\rho(\cdot)$ is Lipschitz, it equals
\begin{align*}
S_k&=\frac{1}{q'd_0}\cdot\frac{1}{k!}\int_{\substack{x^{a}\leq x_i\leq x\\x_1\cdots x_k\leq x}}\frac{\rho\left(\frac{\log \frac{x}{x_1\cdots x_k}}{a\log x}\right)}{x_1\cdots x_k(\log x_1)\cdots (\log x_k)}\,d\mathbf{x}+o_{q'}(1)\\
&=\frac{1}{q'd_0}\cdot\frac{1}{k!}\int_{\substack{a\leq u_1,\ldots, u_k\leq 1\\u_1+\cdots +u_k\leq 1}}\frac{\rho\left(\frac{1-u_1-\cdots -u_k}{a}\right)}{u_1\cdots u_k}d\mathbf{u}+o_{q'}(1),
\end{align*}
where the last integral comes from a change of variables $u_i=\frac{\log x_i}{\log x}$. Combining \eqref{eqq66} with the previous equation, we have shown that
\begin{align*}
\frac{1}{x}\sum_{\substack{x\leq n\leq 2x\\n\equiv c\pmod{q}}}1_{\omega_{>x^{a}}(n)=k}=\frac{1+o_q(1)}{q}\cdot \frac{I_{a,k}}{k!},\,\,\, \frac{1}{x}\sum_{\substack{x\leq n\leq 2x\\n\equiv c\pmod{q}}}1_{\omega_{>x^{b}}(n)=\ell}=\frac{1+o_q(1)}{q}\cdot\frac{I_{b,\ell}}{\ell!},    
\end{align*}
where
\begin{align}\label{eqq99}
I_{\alpha,m}:=\int_{\substack{\alpha\leq u_1,\ldots, u_m\leq 1\\u_1+\cdots +u_m\leq 1}}\frac{\rho\left(\frac{1-u_1-\cdots-u_m}{\alpha}\right)}{u_1\cdots u_m}d\mathbf{u}.
\end{align}
This implies that $g_j\in \mathcal{U}(x,\varepsilon^{-1},\varepsilon)$ for all $x\geq x_0(\varepsilon)$.\\

Now that we have shown that $g_1$ and $g_2$ satisfy our uniform distribution in arithmetic progressions assumption, Theorem \ref{theo_bincorr} with $\omega(X)=\log(3X)$ gives
\begin{align}\label{eqq55a}\begin{split}
\frac{1}{\log_2 x}\sum_{\frac{x}{\log x}\leq n\leq x}\frac{g_1(n)g_2(n+1)}{n}&=\frac{1}{\log_2 x}\sum_{\frac{x}{\log x}\leq n\leq x}\frac{z^{\omega_{>x^{a}}(n)}w^{\omega_{>x^{b}}(n+1)}}{n}\\
&=\left(\frac{1}{x}\sum_{x\leq n\leq 2x}z^{\omega_{>x^{a}}(n)}\right)\left(\frac{1}{x}\sum_{x\leq n\leq 2x}w^{\omega_{>x^{b}}(n)}\right)+o(1).
\end{split}
\end{align}
Note that the numbers $n\in [\frac{x}{\log x},x]$ with $\omega_{>x^{a}}(n)\neq \omega_{>n^{a}}(n)$ have a prime divisor on the interval $[(\frac{x}{\log x})^{a},x^{a}]$, so their contribution to the left-hand side of the above sum is bounded by
\begin{align*}
\sum_{(\frac{x}{\log x})^a\leq p\leq x^{a}}\frac{1}{\log_2 x}\sum_{\substack{\frac{x}{\log x}\leq n\leq x\\p\mid n}}\frac{1}{n}=o(1).
\end{align*}
We can do a similar computation to exclude the terms with $\omega_{>x^{a}}(n)\neq \omega_{>n^{a}}(n)$ on the right-hand side of \eqref{eqq55a}. Applying the same arguments also to $\omega_{>x^{b}}(n)$, \eqref{eqq55a} takes the form
\begin{align}\label{eqq55}
\frac{1}{\log_2 x}\sum_{\frac{x}{\log x}\leq n\leq x}\frac{z^{\omega_{>n^{a}}(n)}w^{\omega_{>n^{b}}(n+1)}}{n}=\left(\frac{1}{x}\sum_{x\leq n\leq 2x}z^{\omega_{>n^{a}}(n)}\right)\left(\frac{1}{x}\sum_{x\leq n\leq 2x}w^{\omega_{>n^{b}}(n)}\right)+o(1).
\end{align}
By the preceding considerations,
\begin{align}\label{eqq89}
\frac{1}{X}\sum_{X\leq n\leq 2X}z^{\omega_{>n^{a}}(n)} =\sum_{0\leq k< \frac{1}{a}}z^k\cdot \frac{I_{a,k}}{k!}+o(1) 
\end{align}
as $X\to \infty$, with $I_{a,k}$ as in \eqref{eqq99}, so summing this dyadically we find that \eqref{eqq89} also holds with the summation range being $1\leq n\leq X$. Thus, by partial summation,
\begin{align}\label{eqq68}\begin{split}
\frac{1}{\log x}\sum_{n\leq x}\frac{z^{\omega_{>n^{a}}(n)}}{n}=\sum_{0\leq k< \frac{1}{a}}z^k\frac{I_{a,k}+o(1)}{k!},\quad\frac{1}{\log x}\sum_{n\leq x}\frac{w^{\omega_{>n^{b}}(n)}}{n}=\sum_{0\leq \ell< \frac{1}{b}}w^{\ell}\frac{I_{b,\ell}+o(1)}{\ell!}.
\end{split}
\end{align}

Based on \eqref{eqq55} and \eqref{eqq68}, if we put
\begin{align*}
c_{k,\ell}(x):&=\frac{1}{\log_2 x}\sum_{\frac{x}{\log x}\leq n\leq x}\frac{1_{\omega_{>n^{a}}(n)=k}1_{\omega_{>n^{b}}(n+1)=\ell}}{n},\\a_k(x):&=\frac{1}{\log x}\sum_{n\leq x}\frac{1_{\omega_{>n^{a}}(n)=k}}{n},\quad b_{\ell}(x):=\frac{1}{\log x}\sum_{n\leq x}\frac{1_{\omega_{>n^{b}}(n)=\ell}}{n},
\end{align*}
then we have
\begin{align*}
\sum_{\substack{0\leq k< \frac{1}{a}\\0\leq \ell< \frac{1}{b}}}c_{k,\ell}(x)z^k w^{\ell}=\left(\sum_{0\leq k< \frac{1}{a}}a_k(x)z^k\right)\left(\sum_{0\leq \ell< \frac{1}{b}}b_{\ell}(x)w^{\ell}\right)+o(1)
\end{align*}
for all $z,w\in [-1,1]$. Expanding out, we see that
\begin{align}\label{eqq56}
\sum_{\substack{0\leq k< \frac{1}{a}\\0\leq \ell< \frac{1}{b}}}(c_{k,\ell}(x)-a_k(x)b_{\ell}(x))z^k w^{\ell}=o(1).
\end{align}
We will show that $c_{k,\ell}(x)=a_k(x)b_{\ell}(x)+o(1)$. Suppose for the sake of contradiction that this is not the case. Then, by compactness, we can find a sequence $x_i$ tending to infinity such that the numbers $D_{k,\ell}:=\lim_{i\to \infty}(c_{k,\ell}(x_i)-a_k(x_i)b_{\ell}(x_i))$ exist, and at least one of them is nonzero. Taking limits in \eqref{eqq56}, we infer
\begin{align*}
\sum_{\substack{0\leq k< \frac{1}{a}\\0\leq \ell< \frac{1}{b}}}D_{k,\ell}z^k w^{\ell}=0
\end{align*}
for all $z,w\in [-1,1]$. We now have a polynomial in two variables vanishing in an open set, so its coefficients $D_{k,\ell}$ must all be zero, which is a contradiction. Thus we have
\begin{align}\label{eqq57}
c_{k,\ell}(x)=a_k(x)b_{\ell}(x)+o(1)=\delta_1^{*}\delta_2^{*}+o(1)
\end{align}
for all $0\leq k< \frac{1}{a}, 0\leq \ell< \frac{1}{b}$, with $\delta_1^{*}:=\frac{I_{a,k}}{k!}$ and $\delta_2^{*}:=\frac{I_{b,\ell}}{\ell!}$.\\

Using \eqref{eqq57} for $x\in \{y_1,y_2,\ldots, y_{J-1}\}$,  where $y_1=x$, $y_{j+1}=\frac{y_j}{\log y_{j}}$ and $y_J\in [\sqrt{\log x},\log x]$, it follows that
\begin{align*}
\frac{1}{\log x}\sum_{n\leq x}\frac{1_{\omega_{>n^{a}}(n)=k}1_{\omega_{>n^{b}}(n+1)=\ell}}{n}&=\frac{1}{\log x}\sum_{j=1}^{J-1}\log \log y_j\cdot (\delta_1^{*}\delta_2^{*}+o(1))\\
&=(\delta_1^{*}\delta_2^{*}+o(1))\frac{1}{\log x}\sum_{j=1}^{J-1}\log \frac{y_j}{y_{j+1}}\\
&=\delta_1^{*}\delta_2^{*}+o(1)=a_k(x)b_{\ell}(x)+o(1)
\end{align*}
by telescopic summation. Taking limits as $x\to \infty$ from this, we reach the statement of the theorem about logarithmic densities.\\

For the part of the theorem involving asymptotic density, we apply the same argument as above, but with $1\leq \omega(X)\leq \log(3X)$ an arbitrary function tending to infinity (instead of $\omega(X)=\log(3X)$). We again have
\begin{align}\label{eqq79}
\frac{1}{\log \omega(x)}\sum_{\frac{x}{\omega(x)}\leq n\leq x}\frac{1_{\omega_{>n^{a}}(n)=k}1_{\omega_{>n^{b}}(n+1)=\ell}}{n}=\delta_1^{*}\delta_2^{*}+o(1).
\end{align}
In particular, we get
\begin{align}\label{eqq59}
\frac{1}{x}\sum_{n\leq x}1_{\omega_{>n^{a}}(n)=k}1_{\omega_{>n^{b}}(n+1)=\ell}\geq \frac{1}{2}\delta_1^{*}\delta_2^{*}\frac{\log \omega(x)}{\omega(x)}
\end{align}
for all large enough $x$ (where large enough depends on the function $\omega(X)$). Now, supposing that the part of Theorem \ref{theo_omega} concerning asymptotic density fails, there is a function $\psi(x)$ tending to infinity such that the left-hand side of \eqref{eqq59} is $\leq \frac{1}{\psi(x)}$ for infinitely many integers $x$. However, taking $\omega(x)=\psi(x)$ in \eqref{eqq59}, we get a contradiction as $x\to \infty$. Hence, there exists some $c_0>0$ such that the left-hand side of \eqref{eqq59} is $\geq c_0$ for all large enough $x$, which was to be shown. \qedd

Our theorems on smooth numbers follow rather quickly from Theorem \ref{theo_omega}. In fact, one could also deduce these applications directly from Theorem \ref{theo_bincorr}, using the fact that smooth numbers are uniformly distributed in arithmetic progressions. We leave the details of this alternative argument to the interested reader.\\

\textbf{Proof of Theorem \ref{theo_density}.} It follows from \eqref{eqq90} with $q'=1$ and partial summation that the set  $\{n\in \mathbb{N}:\,\, P^{+}(n)\leq n^{a}\}$ has logarithmic density $\rho(\frac{1}{a})$. Taking $k=\ell=0$ in Theorem \ref{theo_omega} and noticing that $\omega_{>y}(n)=0$ if and only if $P^{+}(n)\leq y$, the conclusion is immediate.\qedd\\

Theorem \ref{theo_erdos} is a corollary to Theorem \ref{theo_density}, as we will see next.\\

\textbf{Proof of Theorems \ref{theo_erdos} and \ref{theo_integral}.} As mentioned in the introduction, taking $\alpha=0$ in Theorem \ref{theo_integral} implies Theorem \ref{theo_erdos}, since by symmetry 
\begin{align*}
\int_{\substack{(x,y)\in [0,1]^2\\x\geq y}}u(x)u(y)\, dx\, dy=\frac{1}{2}\int_{\substack{(x,y)\in [0,1]^2}}u(x)u(y)\, dx\, dy=\frac{1}{2},    
\end{align*}
where the last equality comes from the fundamental theorem of calculus and the fact that $u(x)=\frac{d}{dx}\rho\left(\frac{1}{x}\right)$. Thus it suffices to prove Theorem \ref{theo_integral}. Let $0<a,b,c,d<1$ be given real numbers with $a<c$ and $b<d$. Applying the inclusion-exclusion formula to the sets $\{n\in \mathbb{N}:\,\, P^{+}(n)\leq n^{a}\},\ldots, \{n\in \mathbb{N}:\,\, P^{+}(n)\leq n^{d}\}$ and employing Theorem \ref{theo_density} and the fundamental theorem of calculus, we see that
\begin{align*}
&\delta(\{n\in \mathbb{N}:\,\, n^{a}<P^{+}(n)<n^{b},\,\, n^{c}<P^{+}(n+1)<n^{d}\})\\
&=\rho\left(\frac{1}{b}\right)\rho\left(\frac{1}{d}\right)-\rho\left(\frac{1}{a}\right)\rho\left(\frac{1}{d}\right)-\rho\left(\frac{1}{b}\right)\rho\left(\frac{1}{c}\right)+\rho\left(\frac{1}{a}\right)\rho\left(\frac{1}{c}\right)\\
&=\left(\rho\left(\frac{1}{d}\right)-\rho\left(\frac{1}{c}\right)\right)\left(\rho\left(\frac{1}{b}\right)-\rho\left(\frac{1}{a}\right)\right)\\
&=\int_{a}^{b}\int_{c}^{d}u(x)u(y)\,dx\,dy.
\end{align*}
In other words, for any rectangle $\mathcal{R}\subset [0,1]^2$ parallel to the coordinate axes we have
\begin{align}\label{eqq100}
\delta\left(\left\{n\in \mathbb{N}:\,\, \left(\frac{\log P^{+}(n)}{\log n},\frac{\log P^{+}(n+1)}{\log n}\right)\in \mathcal{R}\right\}\right)=\int_{\mathcal{R}}u(x)u(y)\,dx\,dy.    
\end{align}
Now, if $S\subset [0,1]^2$ is any set such that $1_S$ is Riemann integrable, we can approximate $S$ from the inside and outside with finite unions of rectangles, so by the monotone convergence theorem we see that \eqref{eqq100} continues to hold for such sets $S$. Taking $S=T_{\alpha}$, Theorem \ref{theo_integral} is proved.\qedd\\

\textbf{Proof of Theorem \ref{theo_hildebrand}.} Let $1\leq \omega(X)\leq \log(3X)$ be a function tending to infinity. Defining $F_u(n):=1_{P^{+}(n)\leq n^{u}}$, by the inclusion-exclusion principle we have
\begin{align*}
&\frac{1}{\log \omega(x)}\sum_{\frac{x}{\omega(x)}\leq n\leq x}\frac{1_{P^{+}(n)\in [n^a,n^b]}1_{P^{+}(n+1)\in [n^c,n^{d}]}}{n}\\
&=\frac{1}{\log\omega(x)}\sum_{\frac{x}{\omega(x)}\leq n\leq x} \frac{F_{b}(n)F_d(n+1)-F_a(n)F_d(n+1)-F_b(n)F_c(n+1)+F_a(n)F_c(n+1)}{n}. \end{align*}
From \eqref{eqq79} (with $k=\ell=0$), it follows that the previous expression is
\begin{align}\label{eqq81}
&\rho\left(\frac{1}{b}\right)\rho\left(\frac{1}{d}\right)-\rho\left(\frac{1}{a}\right)\rho\left(\frac{1}{d}\right)-\rho\left(\frac{1}{b}\right)\rho\left(\frac{1}{c}\right)+\rho\left(\frac{1}{a}\right)\rho\left(\frac{1}{c}\right)+o(1)\nonumber\\
&=\left(\rho\left(\frac{1}{d}\right)-\rho\left(\frac{1}{c}\right)\right)\left(\rho\left(\frac{1}{b}\right)-\rho\left(\frac{1}{a}\right)\right)+o(1).
\end{align}
In particular, as in \eqref{eqq59}, we get
\begin{align}\label{eqq80}
\sum_{n\leq x}\frac{1_{P^{+}(n)\in [n^a,n^b]}1_{P^{+}(n+1)\in [n^c,n^{d}]}}{n}\geq \frac{1}{2} c_0(a,b,c,d)\frac{\log \omega(x)}{\omega(x)},
\end{align}
where $c_0(a,b,c,d)>0$ is the constant in \eqref{eqq81}, and since $\omega(X)$ was allowed to tend to infinity as slowly as we please, the left-hand side of \eqref{eqq80} is lower-bounded by some positive constant, as asserted.\qedd\\

Lastly, we will deduce our quadratic character sum bound from the main theorem.\\

\textbf{Proof of Theorem \ref{theo_burgess}.} The first part of the theorem will follow directly from Theorem \ref{theo_bincorr}, once we show that for any fixed $a,q\in \mathbb{N}$ we have
\begin{align*}
\sum_{\substack{x\leq n\leq 2x\\n\equiv a\pmod{q}}}\chi_Q(n)=o(x)    
\end{align*}
as $x\to \infty$. Denoting $d_0=(a,q)$, $a'=\frac{a}{d_0}$, $q'=\frac{q}{d_0}$, and using complete multiplicativity, it suffices to show that
\begin{align*}
\sum_{\substack{\frac{x}{d_0}\leq m\leq \frac{2x}{d_0}\\m\equiv a'\pmod{q'}}}\chi_Q(m)=o(x).
\end{align*}
Expanding the congruence condition in terms of Dirichlet characters, we are left with showing that
\begin{align}\label{eqq72}
\sum_{\frac{x}{d_0}\leq m\leq \frac{2x}{d_0}}\chi_Q(m)\psi(m)=o(x) \end{align}
for all Dirichlet characters $\psi\pmod{q'}$. Note that the character $\chi^{*}:=\chi_Q\psi$ has modulus $Q^{*}:=Qq'\leq x^{4-\varepsilon/2}$ if $x$ is large enough. In addition, the character $\chi_Q\psi$ cannot be the principal character, since then $\chi_Q$ would be induced by $\psi$, which has modulus $q'<Q$ (since $Q(x)$ is assumed to tend to infinity with $x$), contradicting the assumption that $\chi_Q$ is primitive. The number $Q^{*}$ is not necessarily cube-free, but we can apply a slight generalization of the Burgess bound from \cite[formula (12.56)]{iwaniec-kowalski} to bound the left-hand side of \eqref{eqq72} with
\begin{align*}
\ll_{r,\varepsilon} \left(\frac{x}{d_0}\right)^{1-\frac{1}{r}}q^{\frac{1}{r}}(Q^{*})^{\frac{r+1}{4r^2}+\varepsilon^2}=o(x)    
\end{align*}
for $r=10\lfloor\varepsilon^{-2}\rfloor$, say. Now the first part of the theorem has been proved.\\

For the proof of \eqref{eqq73}, note that the quantity on the left-hand side of that formula is
\begin{align*}
 &\frac{1}{\log x}\sum_{\substack{n\leq x\\(n(n+1),Q)=1}}\frac{1}{n}\cdot \frac{1-\chi_Q(n)}{2}\cdot \frac{1-\chi_Q(n+1)}{2}\\
 &=\frac{1}{4\log x}\sum_{\substack{n\leq x\\(n(n+1),Q)=1}}\frac{1}{n}-\frac{1}{4\log x}\sum_{n\leq x}\frac{\chi_{Q}(n)\chi_0(n+1)}{n}\\
 &\quad -\frac{1}{4\log x}\sum_{n\leq x}\frac{\chi_{0}(n)\chi_Q(n+1)}{n}+\frac{1}{4\log x}\sum_{n\leq x}\frac{\chi_{Q}(n)\chi_Q(n+1)}{n},
\end{align*}
where $\chi_0$ stands for the principal character $\pmod{Q}$. Here the first term equals the right-hand side of \eqref{eqq73} by elementary sieve theory. The other three terms are seen to be $o(1)$ just as in the first part of the theorem (in order to apply Theorem \ref{theo_bincorr}, it suffices that one of $\chi_Q$ and $\chi_0$ is uniformly distributed in arithmetic progressions).\\

For the last part of the theorem, namely proving \eqref{eqq74}, we apply the same argument as in the second part to show that
\begin{align*}
\frac{1}{\log \omega(x)}\sum_{\frac{x}{\omega(x)}\leq n\leq x} \frac{1_{n,n+1\,\, \textnormal{QNR}\pmod{Q}}}{n}=\frac{1}{4}\prod_{p\mid Q}\left(1-\frac{2}{p}\right)+o(1).   
\end{align*}
Since $\omega(X)$ is any function tending to infinity slowly, we can apply exactly the same argument as at the end of the proof of Theorem \ref{theo_omega} to conclude that \eqref{eqq74} holds.\qedd

\appendix

\section{Appendix: Stability of mean values of multiplicative functions}
\label{a:a}

We prove Lemma \ref{le_stability}, which was used in the proof of Theorem \ref{theo_bincorr} and tells that mean values of the functions $g_j$ over the arithmetic progression $a\pmod{q}$ vary very slowly in terms of the interval over which the mean value is taken. The case $q=1$ of the lemma was proved by Elliott \cite{elliott-lipschitz} and refined by Granville and Soundararajan \cite[Proposition 4.1]{gs-spectrum} (see also \cite[Theorem 4]{gs-decay}). Also Matthiesen's work \cite{matthiesen-fourier} contains estimates of the type of Lemma \ref{le_stability}, but for the sake of completeness we give a proof here. We have not aimed to optimize the error terms in the lemma.\\

\textbf{Proof of Lemma \ref{le_stability}.} By writing
\begin{align*}
\frac{1}{x}\sum_{\substack{x\leq n\leq 2x\\n\equiv a\pmod{q}}}g(n)=\frac{1}{x}\sum_{\substack{n\leq 2x\\n\equiv a\pmod{q}}}g(n)-\frac{1}{x}\sum_{\substack{n\leq x\\n\equiv a\pmod{q}}}g(n)   
\end{align*}
and the same with $\frac{x}{y}$ in place of $x$, we see that it suffices to show that
\begin{align}\label{eqq75}
\bigg|\frac{1}{x}\sum_{\substack{n\leq x\\n\equiv a\pmod{q}}}g(n)-\frac{1}{x/y}\sum_{\substack{n\leq \frac{x}{y}\\n\equiv a \pmod{q}}}g(n)\bigg|\ll_q (\log x)^{-\frac{1}{400}}.     
\end{align}
for $y\in [1,2\log^{10} x]$. Putting $d_0:=(a,q)$, $a':=\frac{a}{d_0}$ and $q':=\frac{q}{d_0}$, \eqref{eqq75} becomes
\begin{align*}
\bigg|\frac{1}{x}\sum_{\substack{n'\leq \frac{x}{d_0}\\n'\equiv a'\pmod{q'}}}g(d_0n')-\frac{1}{x/y}\sum_{\substack{n'\leq \frac{x}{d_0y}\\n'\equiv a' \pmod{q'}}}g(d_0n')\bigg|\ll_q (\log x)^{-\frac{1}{400}}.
\end{align*}
Making use of the orthogonality of Dirichlet characters and the triangle inequality, it suffices to show that
\begin{align*}
\bigg|\frac{1}{x}\sum_{n'\leq \frac{x}{d_0}}g(d_0n')\chi(n')-\frac{1}{x/y}\sum_{n'\leq \frac{x}{d_0y}}g(d_0n')\chi(n')\bigg|\ll_q (\log x)^{-\frac{1}{400}}.
\end{align*}
for all Dirichlet characters $\chi\pmod{q'}$. Writing $n'=rm$, where $(m,d_0)=1$ and $r\mid d_0^{\infty}$ (meaning that $r\mid d_0^k$ for some $k$), and using the fact that $g(d_0rm)=g(d_0r)g(m)$, the previous bound will follow from 
\begin{align}\label{eqq76}
\sum_{r\mid d_0^{\infty}}\frac{1}{d_0r}\bigg|\frac{d_0r}{x}\sum_{\substack{m\leq \frac{x}{d_0r}\\(m,d_0)=1}}g(m)\chi(m)-\frac{d_0r}{x/y}\sum_{\substack{m\leq \frac{x}{d_0ry}\\(m,d_0)=1}}g(m)\chi(m)\bigg|\ll_q (\log x)^{-\frac{1}{400}}.   
\end{align}
The terms $r>\log x$ can be discarded, since 
\begin{align*}
\sum_{\substack{r\mid d_0^{\infty}\\r>\log x}}\frac{1}{rd_0}\leq(\log x)^{-\frac{1}{2}}\frac{1}{d_0}\prod_{p\mid d_0}\left(1+\frac{1}{p^{\frac{1}{2}}}+\frac{1}{p}+\cdots \right)\ll_{q}(\log x)^{-\frac{1}{2}}, 
\end{align*}
since $d_0\leq q$. Writing $x':=\frac{x}{d_0r}\gg_q \frac{x}{\log x}$ and applying the triangle inequality to \eqref{eqq76}, together with the simple fact that $\sum_{r\mid d_0^{\infty}}\frac{1}{r}\ll_q 1$, it suffices to show that
\begin{align}\label{eqq77}
\bigg|\frac{1}{x'}\sum_{m\leq x'}g(m)\chi(m)\psi_0(m)-\frac{1}{x'/y'}\sum_{m\leq \frac{x'}{y}}g(m)\chi(m)\psi_0(m)\bigg|\ll_q (\log x)^{-\frac{1}{400}}  
\end{align}
for all $\frac{x}{\log x}\ll_q x'\leq x$, $1\leq y'\leq 2\log^{10} x$ and for all characters $\chi\pmod{q}$, with $\psi_0$ the principal character $\pmod{d_0}$. Note that 
\begin{align}\label{eqq78}
\mathbb{D}(g\chi\psi_0,f;x')=\mathbb{D}(g\chi,f;x')-O_q(1)
\end{align}
for any function $f:\mathbb{N}\to \mathbb{D}$, so we may replace $\psi_0$ with $1$ in any computations involving the pretentious distance.\\

Consider the character $\chi \pmod{q}$ for which the left-hand side of \eqref{eqq77} is maximal. If $\chi$ is complex, we may apply an argument of Granville and Soundararajan (see \cite[Lemma C.1]{mrt}) and the assumption that $g$ is real-valued to obtain the pretentious distance bound
\begin{align*}
\sqrt{M}:=\inf_{|t|\leq x}\mathbb{D}\left(g\chi\psi_0,n^{it};\frac{x'}{y'}\right)=\inf_{|t|\leq x}\mathbb{D}(g\chi,n^{it};x)-O_q(\log_3 x)\geq \frac{1}{10}\sqrt{\log \log x}   
\end{align*}
by \eqref{eqq78}, since $q$ is fixed and $x$ is large enough. Thus by Halász's theorem \cite[Chapter III.4]{tenenbaum}, we may bound \eqref{eqq77} by $\ll Me^{-M}\ll (\log x)^{-1/200}$.\\

In the opposite case that $\chi$ is real in \eqref{eqq77}, we appeal to \cite[Proposition 4.1]{gs-spectrum}, provided that $\mathbb{D}(g\chi\psi_0;1;x)\leq \frac{2}{3}\sqrt{\log \log x}$ holds. This gives a bound of $\ll (\log x)^{-\frac{1}{10}}$ for \eqref{eqq77}, so we may assume that  $\mathbb{D}(g\chi\psi_0;1;x)> \frac{2}{3}\sqrt{\log \log x}$. But since $g\chi$ is real-valued, again by \cite[Lemma C.1]{mrt} we have
\begin{align*}
\sqrt{M}:&=\inf_{|t|\leq x}\mathbb{D}\left(g\chi\psi_0,n^{it};\frac{x'}{y'}\right)=\inf_{|t|\leq x}\mathbb{D}(g\chi,n^{it};x)-O_q(\log_3 x)\\
&\geq \frac{1}{10} \mathbb{D}(g\chi,1;x)-O_q(\log_3 x)\geq \frac{1}{16}\sqrt{\log \log x}   
\end{align*}
for all large enough $x$. Now, again applying Halász's theorem, \eqref{eqq77} is bounded by $\ll Me^{-M}\ll (\log x)^{-1/400}$. This shows that \eqref{eqq77} always holds, which proves the lemma.\qedd

\section{Appendix: Short exponential sum bounds for multiplicative functions}
\label{a:b}

We prove the short exponential sum estimates over major and minor arcs that were employed in the proof of Theorem \ref{theo_bincorr} in Section \ref{sec:circle}. The proofs of both lemmas follow the ideas of Matom\"aki, Radziwi\l{}\l{} and Tao \cite{mrt} for estimating short exponential sum bounds weighted by a multiplicative function, but require some small modifications to the arguments.\\

\textbf{Proof of Lemma \ref{le_minor}.} Since $\alpha \in \mathfrak{m}$, we have the trivial estimate 
\begin{align*}
\left|\sum_{y\leq n\leq y+H}e(\alpha n)\right|\ll \frac{1}{\|\alpha\|}\ll \frac{H}{W}\ll  H(\log H)^{-\frac{1}{10}},
\end{align*}
so by the triangle inequality we may assume that $\delta_1=0$ in \eqref{eqq33}. We introduce the same nicely factorable set $\mathcal{S}:=\mathcal{S}_{P_1,Q_1,X_0,X}$ as in \eqref{eqq50}. By a simple sieve estimate \cite[Lemma 2.3]{mrt}, we have
\begin{align}\label{eqq49}
\sum_{n\leq X+H}(1-1_{\mathcal{S}}(n))\ll \frac{\log \log H}{\log H}X.    
\end{align}
Hence, by the triangle inequality,
\begin{align*}
\frac{1}{X}\int_{X}^{2X}\left|\frac{1}{H}\sum_{y\leq n\leq y+H}g_1(n)(1-1_{\mathcal{S}}(n))e(\alpha n)\right|\, dy\ll \frac{\log \log H}{\log H}\ll  (\log H)^{-\frac{1}{10}}.    
\end{align*}
This means that \eqref{eqq33} has been reduced to
\begin{align}\label{eqq34}
\sup_{\alpha \in \mathfrak{m}}\frac{1}{X}\int_{X}^{2X}\left|\frac{1}{H}\sum_{y\leq n\leq y+H}g_1(n)1_{\mathcal{S}}(n)e(\alpha n)\right|\, dy\ll  (\log H)^{-\frac{1}{10}}. 
\end{align}
This estimate would follow directly from \cite[Section 3]{mrt} (with $d=1$ there), if the function $g_1$ was completely multiplicative, but we will show that the argument goes through even without that assumption.\\

Let $\mathcal{S'}$ be the set of those $n\leq X$ that have a prime factor from each of the intervals $[P_j,Q_j]$ (defined in Definition \ref{def3}) for $j\geq 2$. We have the Ramaré identity
\begin{align*}
g_1(n)1_{\mathcal{S}}(n)&= \sum_{\substack{n=mp\in \mathcal{S}\\P_1\leq p\leq Q_1}}\frac{g_1(mp)1_{\mathcal{S}}(mp)}{|\{P_1\leq p_1\leq Q_1:\,\, p_1\mid n\}|}\\
&=\sum_{\substack{n=mp\\p\nmid m\\P_1\leq p\leq Q_1}}\frac{g_1(m)g_1(p)1_{\mathcal{S'}}(m)}{1+|\{P_1\leq p_1\leq Q_1:\,\, p_1\mid m\}|}+O\left(\sum_{P_1\leq p\leq Q_1}1_{p^2\mid n}\right)\\
&=\sum_{\substack{n=mp\\P_1\leq p\leq Q_1}}\frac{g_1(m)g_1(p)1_{\mathcal{S'}}(m)}{1+|\{P_1\leq p_1\leq Q_1:\,\, p_1\mid m\}|}+O\left(\sum_{P_1\leq p\leq Q_1}1_{p^2\mid n}\right).
\end{align*}
By trivial estimation,
\begin{align*}
\sum_{n\leq X+H}\sum_{P_1\leq p\leq Q_1}1_{p^2\mid n}\ll \sum_{p\geq P_1}\frac{X}{p^2}\ll XW^{-200}\ll X(\log H)^{-\frac{1}{10}},    
\end{align*}
so proving \eqref{eqq34} has been reduced to proving
\begin{align*}
\sum_{P_1\leq p\leq Q_1}\sum_{m}\frac{1_{\mathcal{S}'}(m)g_1(m)g_1(p)e(mp\alpha)}{1+|\{P_1\leq p_1\leq Q_1:\,\, p_1\mid m\}|}\int_{\mathbb{R}}\theta(x)1_{x\leq mp\leq x+H}\,dx \ll HX(\log H)^{-\frac{1}{10}}    
\end{align*}
for all measurable functions $|\theta(x)|\leq 1$ supported on $[0,X]$. This is same expression as in \cite[Section 3]{mrt}, so the proof continues from here in an identical manner (since the rest of the argument does not use multiplicativity).\qedd\\

\textbf{Proof of Lemma \ref{le_major}.} We follow the proof of the major arc exponential sum in \cite[Section 4]{mrt}. However, here we need to be a bit more careful when approximating the exponential $e(\alpha n)$ with $e(\frac{an}{q})$, as we do not want to lose a factor of $\frac{W}{q}$ that would come from a partial summation approximation of $e(\alpha n)$.\\

By our assumption  $g_1\in \mathcal{U}_{\omega}(x,\exp_2(\varepsilon^{-2}),2/\exp_2(\varepsilon^{-2}),\delta_1)$ and formula \ref{eqq49}, it suffices to show that
\begin{align}\label{eqq101}
\frac{1}{X}\int_{X}^{2X}\left|\frac{1}{H}\sum_{x\leq n\leq x+H}(g_1(n)1_{\mathcal{S}}(n)-\delta_1')e(\alpha n)\right|\,dx=o_{\varepsilon\to 0}(1),
\end{align}
where the nicely factorable set $\mathcal{S}$ is as in \eqref{eqq50} and
\begin{align*}
\delta_1':=\frac{1}{X}\sum_{X\leq m\leq 2X}g_1(n)1_{\mathcal{S}}(n)    
\end{align*} 
Let  $H':=\frac{H}{W^3}$. By exchanging the order of integration and summation, we have
\begin{align*}
\frac{1}{H}\sum_{x\leq n\leq x+H}a_n=\frac{1}{H}\int_{x}^{x+H}\frac{1}{H'}\sum_{y\leq n\leq y+H'}a_n\,dy+O\left(\frac{H'}{H}\right)
\end{align*}
for any  $a_n\in \mathbb{D}$. Applying this, we see that the left-hand side of \eqref{eqq101} is
\begin{align*}
&= \frac{1}{X}\int_{X}^{2X}\left|\frac{1}{H}\int_{x}^{x+H}\frac{1}{H'}\sum_{y\leq n\leq y+H'}(g_1(n)1_{\mathcal{S}}(n)-\delta_1')e(\alpha n)\, dy\right|\,dx+O\left(\frac{1}{W^2}\right)\nonumber\\
&\ll \frac{1}{HX}\int_{X}^{2X}\int_{x}^{x+H}\left|\frac{1}{H'}\sum_{y\leq n\leq y+H'}(g_1(n)1_{\mathcal{S}}(n)-\delta_1')e\left(\frac{an}{q}\right)\right|\,dy\,dx+O\left(\frac{1}{W^2}\right),
\end{align*}
where we used the fact that any $n\in [y,y+H']$ obeys
\begin{align*}
e(\alpha n)&=e(\alpha y)\,e(\alpha (n-y))\\
&=e(\alpha y)\,e\left(\frac{a}{q}(n-y)\right)+O\left(\frac{1}{W^2}\right)\\
&=e\left(\left(\alpha-\frac{a}{q}\right)y\right)e\left(\frac{an}{q}\right)+O\left(\frac{1}{W^2}\right)    
\end{align*}
by the inequality
\begin{align*}
\left|e(\alpha (n-y))-e\left(\frac{a}{q}(n-y)\right)\right|\leq 2\pi \left|\alpha-\frac{a}{q}\right||n-y|\leq \frac{2\pi W}{qH}\cdot H'\leq \frac{2\pi}{W^2}.
\end{align*}
By exchanging the order of integration above, it suffices to show that
\begin{align}\label{eqq105}
\frac{1}{X}\int_{X}^{2X}\left|\frac{1}{H'}\sum_{x\leq n\leq x+H'}(g_1(n)1_{\mathcal{S}}(n)-\delta_1')e\left(\frac{an}{q}\right)\right|\,dx\ll \varepsilon.
\end{align}
for all  $1\leq a\leq q\leq W=\log^{5}H$, and with  $H'=\frac{H}{W^3}$, as before. By splitting into residue classes  $\pmod{q}$, \eqref{eqq105} would follow from
\begin{align}\label{eqq107}
 \frac{1}{X}\int_{X}^{2X}\bigg|\frac{1}{H'}\sum_{\substack{x\leq n\leq x+H'\\n\equiv b\pmod{q}}}g_1(n)1_{\mathcal{S}}(n)-\frac{1}{qX}\sum_{X\leq n\leq 2X}g_1(n)1_{\mathcal{S}}(n)\bigg|\,dx\ll \frac{\varepsilon}{q}    
\end{align}
for all  $1\leq a\leq q\leq W$. Applying the triangle inequality and Lemma \ref{le_mr} (and the fact that $q\leq W$), it suffices to show that
\begin{align}\label{eqq108}
 \bigg|\frac{1}{X}\sum_{\substack{X\leq n\leq 2X\\n\equiv b\pmod{q}}}g_1(n)1_{\mathcal{S}}(n)-\frac{1}{qX}\sum_{X\leq n\leq 2X}g_1(n)1_{\mathcal{S}}(n)\bigg|\ll \frac{\varepsilon}{q}.    
\end{align}
As in \cite[Section 2]{mrt}, the fundamental lemma of sieve theory gives for  $q\leq W$ the estimate 
\begin{align*}
\sum_{\substack{X\leq n\leq 2X\\n\equiv a\pmod{q}}}(1-1_{\mathcal{S}}(n))\ll \frac{X}{q}\cdot \frac{\log \log H}{\log H}\ll \frac{\varepsilon}{q}X.
\end{align*}
Taking this into account on both sides of \eqref{eqq108}, that claim is reduced to
\begin{align*}
\bigg|\frac{1}{X}\sum_{\substack{X\leq n\leq 2X\\n\equiv b\pmod{q}}}g_1(n)-\frac{1}{qX}\sum_{X\leq n\leq 2X}g_1(n)\bigg|\ll \frac{\varepsilon}{q}    
\end{align*}
for  $X\in [\frac{x}{\omega(x)},x]$, which follows from our assumption $g_1\in \mathcal{U}_{\omega}(x,\exp_2(\varepsilon^{-2}),2/\exp_2(\varepsilon^{-2}),\delta_1)$ and the fact that  $q\leq \log^{5} H\leq \exp_2(10\varepsilon^{-1})$. The proof is complete.\qedd

\bibliography{myrefs_smooths}{}
\bibliographystyle{plain}

\medskip \medskip
\textsc{Department of Mathematics and statistics, University of Turku, 20014 Turku, Finland}\\
\textit{Email address:} \textup{\texttt{joni.p.teravainen@utu.fi}}

\end{document}